\documentclass[11pt, reqno]{amsart}

 \usepackage{amsmath}
 \usepackage{amssymb}
\usepackage{bm}

\DeclareMathAccent{\mathring}{\mathalpha}{operators}{"17}

 \usepackage{color}

\newcommand{\mysection}[1]{\section{#1}
      \setcounter{equation}{0}}

\newtheorem{theorem}{Theorem}[section]
\newtheorem{lemma}[theorem]{Lemma}

\theoremstyle{definition}
\newtheorem{assumption}{Assumption}[section]

\theoremstyle{remark}
\newtheorem{remark}{Remark}[section]

\newcommand\cbrk{\text{$]$\kern-.15em$]$}}
\newcommand\opar{\text{\raise.2ex\hbox{${\scriptstyle | }$}\kern-.34em$($} }

 \makeatletter
 \def\dashint{%
 \operatorname%
 {\,\,\text{\bf--}\kern-.98em\DOTSI\intop\ilimits@\!\!}}

 \makeatother

\newcommand\bbeta{\text{\raise-.2ex\hbox{$\bm{\beta}$}}}

\newcommand\bR{\mathbb{R}}

\newcommand\bB{\mathbb{B}}

\newcommand\cF{\mathcal{F}}

\newcommand\frA{\mathfrak{A}}
\newcommand\frB{\mathfrak{B}}

\newcommand\frP{\mathfrak{P}}

\newcommand\dist{{\rm dist}\,}

\renewcommand\v{\check{\!\!\!\phantom{x}}\,}

\newcommand\infsup{\operatornamewithlimits{inf\,\,\,sup}}
\newcommand\supinf{\operatornamewithlimits{sup\,\,\,inf}}

\begin{document}

\title[Dynamic programming principle]
{On the dynamic programming principle for
uniformly nondegenerate stochastic differential
games in domains and the Isaacs equations}

\author{N.V. Krylov}
\thanks{The  author was partially supported by
 NSF Grant DNS-1160569}
\email{krylov@math.umn.edu}
\address{127 Vincent Hall, University of Minnesota,
 Minneapolis, MN, 55455}



\begin{abstract}
We prove the dynamic programming principle
for uniformly nondegenerate stochastic differential games
in the framework of time-homogeneous
 diffusion processes considered up to the first exit
time from a domain.
In contrast with previous results established
for constant stopping times we allow arbitrary
stopping times and randomized ones as well.
There is no  assumption about solvability
of the the Isaacs equation in any sense (classical or viscosity).
The zeroth-order ``coefficient'' and the ``free'' term 
are only assumed to be measurable in the space variable.

We also prove that  value functions are uniquely determined
by the functions defining the corresponding Isaacs
equations and thus stochastic games
with the same Isaacs equation have the same value functions. 
\end{abstract}

\maketitle
\mysection{Introduction}

The dynamic programming principle is one of   basic tools
in the theory of controlled diffusion processes. 
It seems to the author that 
Fleming and  Souganidis in \cite{FS89} were the first authors
who proved the dynamic programming principle  with 
nonrandom stopping times
for stochastic differential {\em games\/} in the whole space
on a finite time horizon.
 They used rather involved technical constructions to overcome
some measure-theoretic
difficulties,  a technique somewhat resembling
the one in Nisio \cite{Ni88}, and the theory of viscosity
solutions.

In \cite{Ko09}  Kovats considers time-homogeneous
stochastic differential games
in a ``weak'' formulation
in smooth {\em domains\/} and proves the dynamic programming principle
again with nonrandom stopping times. He uses approximations
of policies by piece-wise constant ones and proceeds
similarly to \cite{Ni88}.
 
\'Swi{\c e}ch in \cite{Sw96} reverses the arguments
in \cite{FS89} and proves the dynamic programming principle
for time-homogeneous stochastic differential games in the whole space
with constant stopping times ``directly''
from knowing  that the viscosity solutions
exist. His method is quite similar to the
so-called verification principle from the theory
of controlled diffusion processes.

It is also worth mentioning the paper \cite{BL08} by Buckdahn and Li
where the dynamic programming principle
for constant stopping times in the time-inhomogeneous setting
in the whole space is derived by using
the theory of backward--forward stochastic equations.

In this paper
we will be only dealing with the dynamic programming principle
for stochastic differential  games and its relation to the 
corresponding Isaacs equations.
 Concerning all other aspect of the theory
of stochastic differential games we refer the reader to 
\cite{BL08}, \cite{FS89}, \cite{Ko09}, \cite{Ni88},
 and  \cite{Sw96}, and the references therein.

In \cite{Kr_1} we adopted the strategy of    \'Swi{\c e}ch
(\cite{Sw96}) which is based on using the fact that in many cases
the Isaacs equation has a sufficiently regular solution.
In \cite{Sw96} viscosity solutions are used and we relied on
classical ones. In the present article no assumptions are
made on the solvability of Isaacs equations. Here we use
a very general result of \cite{Kr12.2}
(see Theorem 1.1 there) about solvability
of {\em approximating\/}
 Isaacs equations and our Theorem \ref{theorem 1.20.1}
implying that the solutions of approximating equations
approximate the value function in the original problem.
Then we basically pass to the limit in the formulas 
obtained in \cite{Kr_1}.

The main emphasis of \cite{FS89}, \cite{Ko09}, \cite{Ni88},
 and \cite{Sw96} is on proving that
the value functions for stochastic differential games are viscosity
solutions of the corresponding Isaacs equations
and the dynamic programming principle is
used just as a tool to do that. In our setting 
the zeroth-order coefficient and the running payoff
function  
can be just measurable and in this situation
neither our methods nor the methods based on the notion of viscosity 
solution seem to be of much help while characterizing 
the value function as a viscosity solution.

Our main future goal is to develop some tools which would allow us
in a subsequent article to show that the value
functions are of class $C^{0,1}$, provided that the data
are there, for {\em possibly degenerate\/} stochastic
differential games without assuming that the zeroth-order
coefficient is large  enough negative.
 On the way to achieve this goal
one of the main steps, apart from proving the dynamic
programming principle, consists of   proving certain representation
formulas like the ones in
 Theorems 3.2 and 3.3 of \cite{Kr_1} in which the process
is not assumed to be uniformly nondegenerate.
Next important ingredient consists of approximations
results stated as Theorem \ref{theorem 1.20.1}
again for possibly degenerated processes.
By combining Theorem 1.1 of \cite{Kr12.2} with
Theorems 3.2 and 3.3 of \cite{Kr_1} and 
 \ref{theorem 1.20.1},  
we then come to one of the main results of the present article,
Theorem \ref{theorem 1.14.1}, about the dynamic
programming principle in a very general form
including stopping and randomized stopping times.

In Theorem \ref{theorem 2.18.1} we assert the H\"older
continuity of the value function in our case
where the zeroth-order coefficient and the 
running payoff
function can be discontinuous.

Theorem \ref{theorem 1.14.1} concerns
time-homogeneous stochastic differential
 games unlike the time inhomogeneous
 in \cite{FS89}
and generalizes the {\em
corresponding\/} results of \cite{Sw96} and \cite{Ko09},
where however  degenerate case is not excluded.

Our Theorem \ref{theorem 2.19.1} shows that the value function
is uniquely defined by the corresponding Isaacs equation
and is independent of the way the equation is represented
as $\sup\inf$ of linear operators (provided that they
satisfy our basic assumptions). This fact in a somewhat
more restricted situation is also noted in Remark 2.4 of 
\cite{Sw96}.

The article is organized as follows.
In Section \ref{section 2.26.3} we state our main results
to which actually, as we pointed out   implicitly above,
belongs Theorem \ref{theorem 1.20.1}. In Section
\ref{section 5.31.1}
 we give a version of Theorem \ref{theorem 1.14.1}
for the whole space.
Then in Section \ref{section 2.27.1} we prove
a very simple result allowing us to compare
the value functions corresponding to different
data.

Sections \ref{section 2.25.1} and \ref{section 3.18.1}
are devoted to deriving
  approximation results. In Section \ref{section 2.25.1}
we consider the approximations from above whereas
in Section \ref{section 3.18.1} from below. The point is that
we know from \cite{Kr12.2} that one can slightly modify
the underlying Isaacs equation in such a way that the modified
equation would have rather smooth solutions.
These  smooth solutions are shown to coincide
with the corresponding value functions, which in addition
satisfy the dynamic programming principle, and the goal
of Sections \ref{section 2.25.1} and \ref{section 3.18.1} is to show
 that when the modification
``fades away'' we obtain the dynamic programming principle
for the original value function. Theorem \ref{theorem 1.20.1}
is proved for the case that the process can degenerate.
Its version for the uniformly
nondegenerate case is given in Section \ref{section 4.14.1}
where we also prove
Theorem \ref{theorem 2.19.1} about the characterization
of the value function by the Isaacs equation. 
In the final short Section \ref{section 4.16.2} we combine
previous results and prove Theorems  \ref{theorem 1.14.1}
and \ref{theorem 2.18.1}.

\mysection{Main results for bounded domains}
                                           \label{section 2.26.3}

Let $\bR^{d}=\{x=(x_{1},...,x_{d})\}$
be a $d$-dimensional Euclidean space and  let $d_{1}\geq d$
and $k\geq1$ be  integers.
Assume that we are given separable metric spaces
  $A$ and $B$  and let,
for each $\alpha\in A$ and $\beta\in B$ 
  the following 
  functions on $\bR^{d}$ be given: 

(i) $d\times d_{1}$
matrix-valued $\sigma^{\alpha\beta}( x)
 =
(\sigma^{\alpha\beta}_{ij}( x))$,

(ii)
$\bR^{d}$-valued $b^{\alpha\beta}( x)=
(b^{\alpha\beta}_{i }( x))$, and

(iii)
real-valued  functions 
$c^{\alpha\beta}( x) $,   
  $f^{\alpha\beta}( x) $, and  
$g(x)$. 

Take a $\zeta\in C^{\infty}_{0}(\bR^{d})$ with unit
integral and for $\varepsilon>0$ introduce
$\zeta_{\varepsilon}(x)=\varepsilon^{-d}\zeta(x/\varepsilon)$.
For locally summable functions $u=u(x)$ on $\bR^{d}$ define
$$
u^{(\varepsilon)}(x)=u*\zeta_{\varepsilon}(x).
$$

\begin{assumption}
                                    \label{assumption 1.9.1}
(i) a) All the above functions are continuous with respect to
$\beta\in B$ for each $(\alpha,x)$ and continuous with respect
to $\alpha\in A$ uniformly with respect to $\beta\in B$
for each $x$. b) These functions are Borel measurable
functions of  $(\alpha,\beta,x)$, the function
$g(x)$ is bounded and uniformly continuous on $\bR^{d}$, and
$c^{\alpha\beta}\geq0$.

(ii) For any $x \in\bR^{d}$  
\begin{equation}
                                              \label{3.27.4}
\sup_{(\alpha,\beta \in A\times B }(
| c^{\alpha\beta}|+|  f^{\alpha\beta}|)( x)<\infty
\end{equation}
and for any $x,y\in\bR^{d}$ 
 and $(\alpha,\beta)\in A\times B $
$$
\|\sigma^{\alpha\beta}( x)-\sigma^{\alpha\beta}( y)\|\leq K_{1}|x-y|,
\quad
  |b^{\alpha\beta}(  x)-b^{\alpha\beta}( y) |\leq K_{1}|x-y|,
$$
$$
  \|\sigma^{\alpha\beta}( x )\|,|b^{\alpha\beta}( x )|
\leq K_{0}
$$
where   $K_{0}$  and $K_{1}$ are
some fixed constants.

(iii) For any bounded domain $D\subset\bR^{d}$ we have
$$
\|\sup_{(\alpha,\beta )\in A\times B }|  f^{\alpha\beta}  
 |\,\|_{L_{d}(D)}
+\|\sup_{(\alpha,\beta )\in A\times B }  c^{\alpha\beta } 
\,\|_{L_{d}(D)}<\infty,
$$
$$
\|\sup_{(\alpha,\beta )\in A\times B }|  f^{\alpha\beta} -
( f^{\alpha\beta})^{(\varepsilon)}|\,\|_{L_{d}(D)}\to0,
$$
$$
 \|\sup_{(\alpha,\beta )\in A\times B }|  c^{\alpha\beta }
-(  c^{\alpha\beta})^{(\varepsilon)}|\,\|_{L_{d}(D)}\to0,
$$
  as 
$\varepsilon\downarrow0$.

(iv) There is a constant $\delta\in(0,1]$
such that for $\alpha\in A$, $\beta\in B$,  
and $x,\lambda\in\bR^{d}$ we have
$$
\delta|\lambda|^{2}\leq a^{\alpha\beta}_{ij}( x)\lambda_{i}
\lambda_{j}\leq \delta^{-1}|\lambda|^{2}.
$$

The reader understands, of course, that the summation
convention is adopted throughout the article.
\end{assumption}

Let $(\Omega,\cF,P)$ be a complete probability space,
let $\{\cF_{t},t\geq0\}$ be an increasing filtration  
of $\sigma$-fields $\cF_{t}\subset \cF_{t}$ such that
each $\cF_{t}$ is complete with respect to $\cF,P$, and let
$w_{t},t\geq0$, be a standard $d_{1}$-dimensional Wiener process
given on $\Omega$ such that $w_{t}$ is a Wiener process
relative to the filtration $\{\cF_{t},t\geq0\}$.

The set of progressively measurable $A$-valued
processes $\alpha_{t}=\alpha_{t}(\omega)$ is denoted by $\frA$. 
Similarly we define $\frB$
as the set of $B$-valued  progressively measurable functions.
By  $ \bB $ we denote
the set of $\frB$-valued functions 
$ \bbeta(\alpha_{\cdot})$ on $\frA$
such that, for any $T\in(0,\infty)$ and any $\alpha^{1}_{\cdot},
\alpha^{2}_{\cdot}\in\frA$ satisfying
\begin{equation}
                                                  \label{4.5.4}
P(  \alpha^{1}_{t}=\alpha^{2}_{t} 
 \quad\text{for almost all}\quad t\leq T)=1,
\end{equation}
we have
$$
P(  \bbeta_{t}(\alpha^{1}_{\cdot})=\bbeta_{t}(\alpha^{2}_{\cdot}) 
\quad\text{for almost all}\quad t\leq T)=1.
$$

For $\alpha_{\cdot}\in\frA$, 
$ \beta_{\cdot} 
\in\frB$,  and $x\in\bR^{d}$ define $x^{\alpha_{\cdot} 
\beta_{\cdot} x}_{t} $ as a unique solution of the It\^o
equation
\begin{equation}
                                             \label{5.11.1}
x_{t}=x+\int_{0}^{t}\sigma^{\alpha_{s}
\beta_{s} }(  x_{s})\,dw_{s}
+\int_{0}^{t}b^{\alpha_{s}
\beta_{s} }(  x_{s})\,ds.
\end{equation}

For a sufficiently smooth function $u=u(x)$ introduce
$$
L^{\alpha\beta} u( x)=a^{\alpha\beta}_{ij}( x)D_{ij}u(x)+
b ^{\alpha\beta}_{i }( x)D_{i}u(x)-c^{\alpha\beta} ( x)u(x),
$$
where, naturally, $D_{i}=\partial/\partial x_{i}$, $D_{ij}=D_{i}D_{j}$.
 Also set
\begin{equation}
                                                     \label{1.16.1}
H[u](x)=\supinf_{\alpha\in A\,\,\beta\in B}
[L^{\alpha\beta} u(x)+f^{\alpha\beta} (x)].
\end{equation}

Denote
$$
\phi^{\alpha_{\cdot}\beta_{\cdot} x}_{t}
=\int_{0}^{t}c^{\alpha_{s}
\beta_{s} }( x^{\alpha_{\cdot} 
\beta_{\cdot}  x}_{s})\,ds.
$$

Next, fix a bounded domain $D\subset\bR^{d}$,
 define $\tau^{\alpha_{\cdot}\beta_{\cdot} x}$ as the first exit
time of $x^{\alpha_{\cdot} 
\beta_{\cdot} x}_{t}$ from $D$, and introduce
\begin{equation}
                                                    \label{2.12.2}
v(x)=\infsup_{\bbeta\in\bB\,\,\alpha_{\cdot}\in\frA}
E_{x}^{\alpha_{\cdot}\bbeta(\alpha_{\cdot})}\big[\int_{0}^{\tau}
f( x_{t})e^{-\phi_{t}}\,dt+g(x_{\tau})e^{-\phi_{\tau}}\big],
\end{equation}
where the indices $\alpha_{\cdot}$, $\bbeta$, and $x$
at the expectation sign are written  to mean that
they should be placed inside the expectation sign
wherever and as appropriate, that is
$$
E_{x}^{\alpha_{\cdot}\beta_{\cdot}}\big[\int_{0}^{\tau}
f( x_{t})e^{-\phi_{t}}\,dt+g(x_{\tau})e^{-\phi_{\tau}}\big]
$$
$$
:=
E \big[
 g(x^{\alpha_{\cdot}\beta_{\cdot}  x}
_{\tau^{\alpha_{\cdot}\beta_{\cdot}  x}}
)
e^{-\phi^{\alpha_{\cdot}\beta_{\cdot}  x}
_{\tau^{\alpha_{\cdot}\beta_{\cdot}  x}}}
+\int_{0}^{\tau^{\alpha_{\cdot}\beta_{\cdot}  x}}
f^{\alpha_{t}\beta_{t}   }
( 
x^{\alpha_{\cdot}\beta_{\cdot}  x}_{t})
e^{-\phi^{\alpha_{\cdot}\beta_{\cdot}  x}_{t}}\,dt\big].
$$
Observe that
this definition makes perfect sense due to 
Theorem 2.2.1  of \cite{Kr77} and
 $v(x)=g(x)$ in $\bR^{d}\setminus D$.

Here is our first main result before which we introduce one more 
assumption. 

\begin{assumption}
                                           \label{assumption 3.19.1}
There exists a nonnegative $G\in C(\bar{D})\cap C^{2}_{loc}(D)$
such that $G=  0$ on $\partial D$ and
$$
  L^{\alpha\beta}G\leq-1
$$
in $D$ for all $\alpha\in A$ and $\beta\in B$.

\end{assumption}

\begin{theorem}
                                            \label{theorem 1.14.1}
 
Under the above assumptions

(i) The function $v(x)$   is bounded and continuous in $\bR^{d}$.

(ii)
Let   $\gamma^{\alpha_{\cdot}\beta_{\cdot}x} $
 be an $\{\cF_{t}\}$-stopping
time defined for each $\alpha_{\cdot}\in\frA$, $\beta_{\cdot}\in\frB$,
and $x\in \bR^{d}$ and such
 that $\gamma^{\alpha_{\cdot}\beta_{\cdot}x}\leq 
\tau^{\alpha_{\cdot}\beta_{\cdot} x}$. Also let $\lambda_{t}^{\alpha_{\cdot}
\beta_{\cdot}x}\geq0$ be progressively measurable  functions on $\Omega
\times[0,\infty)$ 
defined for each $\alpha_{\cdot}\in\frA$, $\beta_{\cdot}\in\frB$,
and $x\in \bR^{d}$ and
such that they have finite integrals over finite time intervals
(for any $\omega$).
Then for any $x$
\begin{equation}
                                                           \label{1.14.1}
v(x)=\infsup_{\bbeta\in\bB\,\,\alpha_{\cdot}\in\frA}
E_{x}^{\alpha_{\cdot}\bbeta(\alpha_{\cdot})}\big[
v(x_{\gamma})e^{-\phi_{\gamma}-\psi_{\gamma}}
+\int_{0}^{\gamma}
\{f( x_{t})+\lambda_{t}v(x_{t})\}e^{-\phi_{t}-\psi_{t}}\,dt \big],
\end{equation}
where inside the expectation sign
$\gamma=\gamma^{\alpha_{\cdot}\bbeta(\alpha_{\cdot})x}
$ and
$$
\psi^{\alpha_{\cdot}\beta_{\cdot} x}_{t}
=\int_{0}^{t}
\lambda^{\alpha_{\cdot}\beta_{\cdot} x}_{s}\,ds.
$$

\end{theorem}

\begin{remark}
The above setting is almost identical to that
of \cite{Kr_1} and statement of
Theorem \ref{theorem 1.14.1} is almost identical to that of Theorem
2.2 of \cite{Kr_1}. However, here we did not impose a quite strong
assumption from \cite{Kr_1} that $D$ be approximated
by domains in which the Isaacs equation has
  regular solutions.  On the other hand, we pay for that
by excluding parameters $p$, which are present
in Theorem
2.2 of \cite{Kr_1} and will reappear in our Theorem 
\ref{theorem 2.19.1}.
\end{remark}

Note that
the possibility to vary $\lambda$ in Theorem \ref{theorem 1.14.1}
might be useful while considering stochastic differential games
with stopping in the spirit of \cite{Kr70}.

\begin{theorem} 
                                             \label{theorem 2.18.1}
The function $v$ is locally H\"older continuous
in $D$ with exponent $\theta\in(0,1)$ depending only
on $d$ and $\delta$.

\end{theorem}

Next, we    state a comparison result, for which we need
some new objects and additional assumptions.
 Take an integer $k \geq d$
and assume that on $\bR^{k }$ we are given a mapping
$$
\Pi:\check x\in\bR^{k}\to\Pi(\check x)\in\bR^{d}
$$
which is twice continuously differentiable with bounded
and uniformly continuous
first- and second-order derivatives.

The reader understands that
the case $k =d$ is not excluded  in which case
  $\Pi(\check{x})\equiv\check{x}$ is allowed.

Assume that we are given a separable metric space 
   $P$  and let,
for each $\alpha\in A$, $\beta\in B$, and $p\in P$,
  the following 
  functions on $\bR^{k}$ be given: 

(i) $k\times d_{1}$
matrix-valued $\check \sigma^{\alpha\beta}(p,\check x)
 =
(\check \sigma^{\alpha\beta}_{ij}(p,\check x))$,

(ii)
$\bR^{k}$-valued $\check b^{\alpha\beta}(p,\check x)=
(\check b^{\alpha\beta}_{i }(p,\check x))$, and

(iii)
real-valued  functions $\check r^{\alpha\beta}(p,\check x) $ ,
$\check c^{\alpha\beta}(p,\check x) $, and 
  $\check f^{\alpha\beta}(p,\check x) $.

As usual we introduce
$$
\check{a}^{\alpha\beta}(p,\check x)=(1/2)
\check\sigma^{\alpha\beta}(p,\check x)
(\check\sigma^{\alpha\beta}(p,\check x))^{*}
$$
and for a fixed $\bar p\in P$  define
$$
(\bar a,\bar \sigma,\bar b,\bar c,\bar f,\bar r)^{\alpha\beta}(\check x)
=(\check a,\check\sigma,\check b, \check c,\check f
,\check r)^{\alpha\beta}(\bar{p},\check x).
$$

\begin{assumption}
                                       \label{assumption 5.14.1}
(i) All the above functions apart from $\check r$ are continuous with respect to
$\beta\in B$ for each $(\alpha,p,\check x)$ and continuous with respect
to $\alpha\in A$ uniformly with respect to $\beta\in B$
for each $(p,\check x)$. Furthermore, they are Borel measurable functions
of $(p,\check x)$ for each
$(\alpha,\beta)$    and
$\check c^{\alpha\beta}\geq0$.

(ii) The functions $\bar{\sigma}^{\alpha\beta}(\check  x )$ and
$\bar{b}^{\alpha\beta}(\check  x )$ are uniformly continuous with
respect to $\check x$ uniformly with respect to 
$(\alpha,\beta )\in A\times B $
and
for any $\check x \in\bR^{k}$ 
 and $(\alpha,\beta,p)\in A\times B\times P$
$$
  \|\sigma^{\alpha\beta}(p,\check x )\|,|b^{\alpha\beta}(p,\check x )|
\leq K_{0}.
$$

(iii) We have $\bar{r}\equiv1$ and 
there is a constant $\check{\delta}_{1}\in(0,1]$
such that on
$A\times B\times 
P\times \bR^{k}$ we have
\begin{equation}
                                               \label{6.14.4}
 \check r^{\alpha\beta}(p,\check x)\in  
[\check\delta_{1},\check\delta_{1}^{-1}],\quad 
\check f ^{\alpha\beta}(p,\check x)
=\check r^{\alpha\beta}(p,\check x)  
 \bar f 
^{\alpha\beta}(\check x) .
\end{equation}

(iv) The functions  $c^{\alpha\beta}(  x)$ and 
$f^{\alpha\beta}(  x)$ are
 bounded on $A\times B\times \bR^{d}$.
(This part bears on the objects introduced before
Theorem \ref{theorem 1.14.1}.)

(v)  For any $\check x \in\bR^{k}$  
\begin{equation}
                                              \label{3.27.04}
\sup_{(\alpha,\beta \in A\times B }(
|\bar c^{\alpha\beta}|+|\bar f^{\alpha\beta}|)(\check x)<\infty.
\end{equation}
 
\end{assumption}

A  
function $p^{\alpha_{\cdot}\beta_{\cdot}}_{t}=
p^{\alpha_{\cdot}\beta_{\cdot}}_{t}(\omega)$ given
on $\frA\times\frB\times\Omega\times(0,\infty)$
is said to be {\em  control adapted\/} if,
for any $(\alpha_{\cdot},\beta_{\cdot})\in\frA\times\frB$
it is progressively measurable in $(\omega,t)$ and, for any
$T\in(0,\infty)$, we have
$$
P( p^{\alpha^{1}_{\cdot}\beta^{1}_{\cdot}}_{t}
=p^{\alpha^{2}_{\cdot}\beta^{2}_{\cdot}}_{t}
\quad\text{for almost all}\quad t\leq T)=1
$$
as long as
$$
P(  \alpha^{1}_{t}=\alpha^{2}_{t} ,
\beta^{1}_{t}=\beta^{2}_{t}
\quad\text{for almost all}\quad t\leq T)=1.
$$
The set of $P$-valued control adapted processes is denoted by $\frP$.

  We discussed a way in which
control adapted processes appear naturally
in Remark   2.3 of \cite{Kr_1}.

Fix a $p\in\frP$ and
 for   $\alpha_{\cdot}\in\frA$, $\beta_{\cdot}\in\frB$,
and $\check{x}\in\bR^{k }$ consider the following equation

\begin{equation}
                                                    \label{5.12.1}
\check{x}_{t}=\check{x}+
\int_{0}^{t}\check{\sigma}^{\alpha_{s}\beta_{s}}
(p^{\alpha_{\cdot}\beta_{\cdot}}_{s},\check{x}_{s})\,dw_{s}+\int_{0}^{t}\check{b}^{\alpha_{s}\beta_{s}}
(p^{\alpha_{\cdot}\beta_{\cdot}}_{s},\check{x}_{s})\,ds.
\end{equation}

\begin{assumption}
                                           \label{assumption 5.23.2}
Equation \eqref{5.12.1} satisfies the usual hypothesis, that is
for any $\alpha_{\cdot}\in\frA$, 
$ \beta_{\cdot} 
\in\frB$, and $\check{x}\in\bR^{k }$ it has a unique 
solution
denoted by $ 
\check{x}_{t}^{\alpha_{\cdot}\beta_{\cdot}\check{x}}$ and
 $ \check{x}_{t}^{\alpha_{\cdot}\beta_{\cdot}\check{x}}$ 
is a control adapted process
for each $\check x$.

\end{assumption}

In order to state additional assumptions,
 we need a  possibly unbounded domain  
$ D^{\v} \subset \bR^{k }$ such that 
$$
\Pi( D^{\v})=D.
$$

Denote by
 $
\check \tau^{\alpha_{\cdot}\beta_{\cdot}\check{x}}
 $
  the first exit time of $\check{x}_{t}^{\alpha_{\cdot}\beta_{\cdot}\check{x}}$
from $D^{\v}$ and set
$$
\check \phi^{\alpha_{\cdot}\beta_{\cdot} x}_{t}
=\int_{0}^{t}\check c^{\alpha_{s}
\beta_{s} }( p^{\alpha_{\cdot}\beta_{\cdot}}_{s},
\check x^{\alpha_{\cdot} 
\beta_{\cdot}  x}_{s})\,ds,
$$

Next, suppose that for each $\varepsilon>0$ we are given 
  real-valued Borel measurable
functions $\bar{c}^{\alpha\beta}_{\varepsilon}(\check{x})$ and
$\bar{f}^{\alpha\beta}_{\varepsilon}(\check{x})$ 
defined on $A\times B\times
\bR^{k}$
 and impose

\begin{assumption}
                                   \label{assumption 3.25.1}
(i) For each $\varepsilon>0$   
the functions 
$(\bar{c}_{\varepsilon},\bar{f}_{\varepsilon})^{\alpha\beta} $ 
are 
bounded on $A\times B\times \bar{D^{\v}}$ and uniformly
 continuous with respect to $\check{x}\in\bar{D^{\v}}$
uniformly with respect to $\alpha,\beta$.

 (ii) For any $\check{x}\in D^{\v}$
$$
\sup_{(\alpha_{\cdot},\beta_{\cdot} )\in\frA\times\frB
 }
E^{\alpha_{\cdot}\beta_{\cdot}}_{\check x}
\int_{0}^{\check \tau}
\sup_{\alpha\in A,\beta\in B}
|\bar  c^{\alpha\beta} -\bar c^{\alpha\beta} _{\varepsilon}|
(\check{x}_{t})
e^{-\check{\phi}_{t}}\,dt\to0,
$$
\begin{equation}
                                              \label{3.27.2}
\sup_{(\alpha_{\cdot},\beta_{\cdot} )\in\frA\times\frB
 }
E^{\alpha_{\cdot}\beta_{\cdot}}_{\check x}
\int_{0}^{\check \tau}
\sup_{\alpha\in A,\beta\in B}
|\bar f^{\alpha\beta} -\bar f^{\alpha\beta} _{\varepsilon}|
(\check{x}_{t})
e^{-\check{\phi}_{t}}\,dt\to0
\end{equation}
with the second convergence in \eqref{3.27.2} being uniform
in $D^{\v}$.

(iii) There
exists a  constant  $\check{\delta}\in(0,1]$ such that for 
$\check{x}\in\bR^{k}$, $p\in P$, $\alpha\in A$,
$\beta\in B$, and $\lambda\in\bR^{d}$ we have
$$
\check{\delta}|\lambda|^{2}\leq 
\big|\lambda^{*}\frac{\partial \Pi}{\partial\check x}(\check x)
\check{\sigma}^{\alpha\beta}
(p,\check{x})\big|^{2}\leq\check{\delta}^{-1}|\lambda|^{2}.
$$

\end{assumption}

\begin{remark}
                                           \label{remark 6.4.1}
Assumption \ref{assumption 3.25.1} (iii) is equivalent to
 saying that for  solutions of
\eqref{5.12.1} the processes $\Pi(\check x_{t})$ 
are uniformly nondegenerate.
\end{remark}
 
It is convenient to always lift functions $u$ given
on $\bR^{d}$ to functions given on $\bR^{k}$ by the formula
\begin{equation}
                                                      \label{3.23.1}
u(\check{x}):=u(\Pi(\check{x})).
\end{equation}

For  sufficiently smooth functions $u=u(\check x)$ introduce
$$
\check L^{\alpha\beta} u(p,\check x)=
\check a^{\alpha\beta}_{ij}(p,\check x)
D_{ij}u(\check x)+
\check b ^{\alpha\beta}_{i }(p,x)D_{i}u(\check x)-
\check c^{\alpha\beta} (p,
\check x)u(\check x),
$$
$$
\bar{L}^{\alpha\beta}u(\check{x})
=\check{L}^{\alpha\beta}u(\bar{p},
\check{x} )
$$
(naturally, $D_{i}=\partial/\partial 
\check x_{i}$, $D_{ij}=D_{i}D_{j}$).
 Also set
$$
\check H[u](\check x)=\supinf_{\alpha\in A\,\,\beta\in B}
[\bar L^{\alpha\beta} u(\check x)+\bar f^{\alpha\beta} (\check x)].
$$

\begin{assumption}
                                        \label{assumption 5.24.1}
There
 exists a nonnegative (bounded) function
$\check{G}\in C(\bar{D^{\v}})\cap C^{2}_{loc}(D^{\v})$
such that $\check{G}(\check{x})\to  0$ as
 $\check{x}\in \bar{D^{\v}}$ and
 $\dist(\Pi (\check{x}),\partial D)\to0$
 ($\check{G} =  0$ on $\partial D$ if $k =d$
and $\Pi(\check x)\equiv\check x$)
and
$$
  \check{L}^{\alpha\beta}\check{G}(p,\check{x})\leq-1
$$
in $P\times D^{\v}$ for all $\alpha\in A$ and $\beta\in B$. 
\end{assumption}

Next, take a real-valued function  
$\psi$ on $\bR^{k}$ with finite
 $C^{2}(\bR^{k})$-norm and   introduce 
$$
\check v(\check x)=
\infsup_{\bbeta\in\bB\,\,\alpha_{\cdot}\in\frA}
E_{\check x}^{\alpha_{\cdot}\bbeta(\alpha_{\cdot})}
\big[\int_{0}^{\check\tau}
\check f( p_{t},\check x_{t})e^{-\check\phi_{t}}\,dt+
\psi(\check{x}_{\check\tau})v(  \check{x}_{\check\tau})
e^{-\check\phi_{\check\tau}}\big],
$$
where, naturally, 
$v$ is taken from Theorem \ref{theorem 1.14.1}. 
Assumption \ref{assumption 3.25.1} (iii)
(and the boundedness of $D$) and
Theorem 2.2.1  of \cite{Kr77} allow  us to conclude that that $
P^{\alpha_{\cdot}\beta_{\cdot}}_{x}
(\check \tau^{\alpha_{\cdot}\beta_{\cdot}x}<\infty)=1$.
 Also notice that \eqref{6.14.4}  and
Assumptions \ref{assumption 5.24.1} imply that
for any $\check x\in D^{\v}$
$$
\delta_{1}\sup_{(\alpha_{\cdot}\beta_{\cdot})\in\frA\times\frB}
E^{\alpha_{\cdot}\beta_{\cdot}}_{\check x}
\int_{0}^{\check \tau}|f (p_{t},
\check x_{t})|e^{-\check\phi_{t}}\,dt
$$
$$
\leq 
\sup_{(\alpha_{\cdot}\beta_{\cdot})\in\frA\times\frB}
E^{\alpha_{\cdot}\beta_{\cdot}}_{\check x}
\int_{0}^{\check \tau}|\bar f (
\check x_{t})|e^{-\check\phi_{t}}\,dt
$$
$$
\leq 
\sup_{(\alpha_{\cdot}\beta_{\cdot})\in\frA\times\frB}
E^{\alpha_{\cdot}\beta_{\cdot}}_{\check x}
\int_{0}^{\check \tau}|\bar f (
\check x_{t})-f_{\varepsilon} (
\check x_{t})|e^{-\check\phi_{t}}\,dt+\check G(\check x)
\sup_{\alpha,\beta,\check y}|f^{\alpha\beta}_{\varepsilon}
(\check y)|,
$$
which is finite at least for small $\varepsilon>0$
owing to \eqref{3.27.2}. Hence,  $\check v$ is well defined.

By the way,
observe also that, if $k=d$ and $\Pi(\check x)\equiv\check x$, then
 $\psi(\check{x})v(  \check{x})  =\psi(x)g(x)$
on $\partial D^{\v}=\partial D$.
 
\begin{assumption}
                                         \label{assumption 5.9.1}
For any function $u\in C^{2}_{loc}(D)$ (not $C^{2}_{loc}(D^{\v})$),
 the function
$\psi(\check{x})u(\check{x})$ is $p$-insensitive
in $D^{\v}$ relative to
 $(\check r^{\alpha\beta}(p,\check{x}),
 \check{L}^{\alpha\beta}(p,\check{x}))$
in the terminology of \cite{Kr_1}, that is,
for any $\alpha_{\cdot}$, $\beta_{\cdot}$, and $\check x$
 we have
$$
d\big[
(\psi u)(\check{x}_{t}^{\alpha_{\cdot}\beta_{\cdot}\check{x}})
e^{-\check{\phi}^{\alpha_{\cdot}\beta_{\cdot} \check x}_{t}}\big]
=dm_{t}
$$
$$
+e^{-\check{\phi}^{\alpha_{\cdot}\beta_{\cdot} \check x}_{t}}
\check{r}^{\alpha_{t}\beta_{t}}(p^{\alpha_{\cdot}\beta_{\cdot}}_{t},
\check{x}_{t}^{\alpha_{\cdot}\beta_{\cdot}\check{x}})
\bar{L}^{\alpha_{t}\beta_{t}}(\psi u)(
\check{x}_{t}^{\alpha_{\cdot}\beta_{\cdot}\check{x}})\,dt,
$$
whenever 
$t<\check \tau^{\alpha_{\cdot}\beta_{\cdot}\check{x}}$,
where $m_{t}$ is a local martingale starting at zero.
\end{assumption}

We discuss this assumption in Remark \ref{remark 5.24.1}.
  
Finally, take some $\{\cF_{t}\}$-stopping
times 
$\gamma^{\alpha_{\cdot}\beta_{\cdot}\check{x}} $
   and 
progressively measurable  functions 
$\lambda_{t}^{\alpha_{\cdot}\beta_{\cdot}\check{x}}\geq0$ 
on $\Omega\times[0,\infty)$ 
defined for each $\alpha_{\cdot}\in\frA$, $\beta_{\cdot}\in\frB$, and $\check{x}
\in\bR^{k}$ and
such that $\lambda_{t}^{\alpha_{\cdot}\beta_{\cdot}\check{x}}$ 
  have finite integrals over finite time intervals
(for any $\omega$).
  Introduce
$$
 \psi ^{\alpha_{\cdot}\beta_{\cdot}\check{x}}_{t}
=\int_{0}^{t}\lambda_{s}^{\alpha_{\cdot}\beta_{\cdot}\check{x}}\,ds.
$$

In the following   theorem by quadratic functions
we mean quadratic functions on $\bR^{d}$ (not $\bR^{k }$) (and if $u$
is a function defined in $D$ then we extend it to a function
in a domain in  $\bR^{k}$ following notation \eqref{3.23.1}).

\begin{theorem}
                                        \label{theorem 2.19.1}

 (i) If
for any   $\check{x} \in D^{\v}$   and 
quadratic function $u$, we have
\begin{equation}
                                            \label{3.19.1}
H[u](\Pi (\check{x}))\leq0\Longrightarrow 
\check{H}[u\psi ](\check{x})\leq0,
\end{equation} 
then $\check{v}\leq \psi v$ in $\bR^{k}$ and
for any $\check x\in \bR^{k}$
$$
v\psi(\check{x})
\geq\infsup_{\bbeta\in\bB\,\,\alpha_{\cdot}\in\frA}
E_{\check{x}}^{\alpha_{\cdot}\bbeta(\alpha_{\cdot})}\big[
v\psi(\check{x}_{\gamma\wedge\check \tau})
 e^{-\check{ \phi}_{\gamma\wedge\check \tau}- \psi _{\gamma\wedge
\check \tau}}
$$
\begin{equation}
                                            \label{4.2.1}
+\int_{0}^{\gamma\wedge\check \tau}
\{\check{f}(p_{t},\check{x} _{t} )+\lambda_{t}
v\psi(\check{x} _{t})\}e^{-\check{ \phi}_{t}- \psi _{t}}\,dt \big].
\end{equation}

 (ii) If
for any   $\check{x} \in D^{\v}$  and 
quadratic function $u$, we have
\begin{equation}
                                            \label{3.19.2}
H[u](\Pi (\check{x}))\geq0\Longrightarrow 
\check{H}[u\psi ](\check{x})\geq0,
\end{equation} 
then $\check{v}\geq \psi v$ in $\bR^{k}$ and
for any $\check x\in \bR^{k}$
$$
v\psi(\check{x})
\leq\infsup_{\bbeta\in\bB\,\,\alpha_{\cdot}\in\frA}
E_{\check{x}}^{\alpha_{\cdot}\bbeta(\alpha_{\cdot})}\big[
v\psi(\check{x}_{\gamma\wedge\check \tau})
 e^{-\check{ \phi}_{\gamma\wedge\check \tau}- \psi _{\gamma\wedge
\check \tau}}
$$
\begin{equation}
                                            \label{4.2.r}
+\int_{0}^{\gamma\wedge\check \tau}
\{\check{f}(p_{t},\check{x} _{t} )+\lambda_{t}
v\psi(\check{x} _{t})\}e^{-\check{ \phi}_{t}- \psi _{t}}\,dt \big].
\end{equation}

\end{theorem}

\begin{remark}
                                             \label{remark 5.29.5}

Under the assumptions of Theorem \ref{theorem 1.14.1}
suppose that $c$ and $f$ are bounded.
Take a global barrier $\Psi$,
which is an infinitely differentiable function
on $\bR^{d}$ such that $\Psi\geq1$ on $\bR^{d}$
 and $(L^{\alpha\beta}+c^{\alpha\beta})
\Psi\leq-1$  on $D$ for all $\alpha,\beta$. The existence
of such functions is a simple and well-known fact.

In Theorem \ref{theorem 2.19.1} take $k=d$, $D^{\v}=D$,
and independent of $p$ functions $\check{r}\equiv1$,
$$
\check{\sigma}^{\alpha\beta}(x)=\Psi^{1/2}(x)
\sigma^{\alpha\beta}(x),\quad\check{b}^{\alpha\beta}(x)
=\Psi(x)b^{\alpha\beta}(x)+2a^{\alpha\beta}(x)D\Psi(x),
$$
$$
\check{c}^{\alpha\beta}(x)=-L^{\alpha\beta}\Psi(x),\quad
\check{f}^{\alpha\beta}(x)= f^{\alpha\beta}(x),
\quad\check{g}(x)=\Psi^{-1}(x)g(x),
$$
where $D\Psi$ is the gradient of $\Psi$ (a column vector).

A simple computation shows that
$$
\check{L}^{\alpha\beta}u(x)+\check{f}^{\alpha\beta}
=L^{\alpha\beta}(u\Psi)(x)+f^{\alpha\beta}(x)
$$
and therefore both conditions in 
\eqref{3.19.1} and \eqref{3.19.2} are satisfied with $\psi=\Psi^{-1}$
and by Theorem \ref{theorem 2.19.1} we conclude that
 $\check v=\Psi^{-1} v$. It is still probably worth noting
that to check Assumption \ref{assumption 3.25.1} in this case
we take $(\bar 
c_{\varepsilon}, \bar 
f_{\varepsilon})^{\alpha\beta} =
[( \check c , \check f)^{\alpha\beta}]^{(\varepsilon)}$

This simple observation sometimes helps introducing a new $c\geq1$
when the initial one was zero.
\end{remark}

\begin{remark}
                                          \label{remark 4.2.1}
If $\check{a},\check{b},\check{c}$, and $\check{f}$
are independent of $p$ and
 $k =d$, $\Pi(x)\equiv x$,  
 and $\psi\equiv1$, then Theorem \ref{theorem 2.19.1}
implies that $v=\check{v}$ whenever the functions $H$ and $\check{H}$
coincide. Therefore, $v$ and $\check v$ are uniquely defined by $H$
and not by its particular representation \eqref{1.16.1}
and, for that matter, not by the choice of probability
space, filtration, and the
 Wiener process including its dimension.
By Theorem \ref{theorem 2.19.1} we also have that $v=\check{v}$ 
if $k =d$, $\Pi(x)\equiv x$, and if
 $\check{a},\check{b},\check{c}$, and $\check{f}$ do depend on $p$
but in such a way that
$$
(\check{a},\check{b},\check{c},\check{f})(p,x)
=\check{r}(p,x)(a,b,c,f)(x)
$$
since in that case any smooth function is $p$-insensitive.
In such a situation we see that $\check{v}$ is independent of
 $p\in\frP$ as well.

Also notice that, if in Theorem \ref{theorem 1.14.1} the functions
 $c$ and $f$ are bounded
(see Assumption \ref{assumption 5.14.1} (iv))
 and one takes $k=d$, assumes that the checked functions
are independent of $p$, and finally takes the 
checked functions equal to the unchecked ones and $(\bar 
c_{\varepsilon}, \bar 
f_{\varepsilon})^{\alpha\beta} =
[(  c , f)^{\alpha\beta}]^{(\varepsilon)}$,
  then one sees that assertion (ii) of Theorem \ref{theorem 1.14.1}
follows immediately from 
Theorem 2.2.1  of \cite{Kr77} and Theorem \ref{theorem 2.19.1}.

\end{remark}
\begin{remark}
                                          \label{remark 6.4.01}
Here we discuss the possibility to use dilations.
Take a constant $\mu>0$ and consider the following
modification of \eqref{5.11.1}
\begin{equation}
                                             \label{6.4.3}
x_{t}=x+\int_{0}^{t} \sigma^{\alpha_{s}
\beta_{s} }( \mu x_{s})\,dw_{s}
+\int_{0}^{t} \mu b^{\alpha_{s}
\beta_{s} }( \mu x_{s})\,ds.
\end{equation}
The solution of this equation is denoted by $x^{\alpha_{\cdot}\beta_{\cdot}
x}_{t}(\mu)$. Then let
$$
\phi^{\alpha_{\cdot}\beta_{\cdot} x}_{t}(\mu)
=\int_{0}^{t}\mu^{2}c^{\alpha_{s}
\beta_{s} }(\mu x^{\alpha_{\cdot} 
\beta_{\cdot}  x}_{s}(\mu))\,ds,
$$
denote by $\tau^{\alpha_{\cdot}\beta_{\cdot} x}(\mu)$
the first exit time of
$x^{\alpha_{\cdot}\beta_{\cdot}
x}_{t}(\mu)$ from $\mu^{-1}D$, and set
$$
v(x,\mu)=\infsup_{\bbeta\in\bB\,\,\alpha_{\cdot}\in\frA}
E_{x}^{\alpha_{\cdot}\bbeta(\alpha_{\cdot})}
\big[\int_{0}^{\tau(\mu)}
\mu^{2}f(\mu x_{t}(\mu))e^{-\phi_{t}(\mu)}\,dt
$$
$$
+
g(\mu x_{\tau}(\mu))e^{-\phi_{\tau(\mu)}(\mu)}\big].
$$

A simple application of Theorem \ref{theorem 2.19.1}
with $\Pi(x)=\mu x$ and $\psi\equiv1$ shows that $v(\mu x)=v(x,\mu)$.
Of course, other types of changing the coordinates
are also covered by Theorem \ref{theorem 2.19.1}.
\end{remark}
\begin{remark}
                                               \label{remark 5.24.1}

The case $k >d$ will play a very important role
in a subsequent article about   stochastic differential games.
To illustrate one of applications consider the one-dimensional 
Wiener process $w_{t}$, define $\tau_{x}$
as the first exit time of $x+ w_{t}$ from $(-1,1)$
and introduce
$$
v(x)=E\int_{0}^{\tau_{x}}f(x+w_{t})\,dt,
$$
so that the corresponding (Isaacs) equation becomes 
$$
H[v]:=(1/2)D^{2}v+f=0
$$
in $(-1,1)$ with zero boundary data at $\pm1$. 
We want to show how Theorem \ref{theorem 2.19.1}
allows one to derive the following
\begin{equation}
                                                \label{3.18.6}
v(x)=E\int_{0}^{\check{\tau}_{x}}e^{-w_{t}-(1/2)t}f(x+w_{t}+t)\,dt,
\end{equation}
where $\check{\tau}_{x}$ is the first exit time of $x+w_{t}+t$
from $(-1,1)$. (Of course, \eqref{3.18.6}
  is a simple corollary of Girsanov's theorem.)

In order to do that consider the two-dimensional diffusion
process given by
\begin{equation}
                                                \label{3.18.5}
dx_{t}=dw_{t}+dt,\quad dy_{t}=-y_{t}\,dw_{t}
\end{equation}
starting at 
$$
(x,y)\in D^{\v}_{\varepsilon}=
(-1,1)\times(\varepsilon,\varepsilon^{-1}),
$$ 
where $\varepsilon\in(0,1)$,
let $\tau^{\varepsilon}_{x,y}$ be the first time the process
exits from this domain, and introduce
$$
\check{v}(x,y)=E\big[\int_{0}^{\tau^{\varepsilon}_{x,y}}y_{t}f(x_{t})\,dt
+y_{\tau^{\varepsilon}_{x,y}}v(x_{\tau^{\varepsilon}_{x,y}})\big].
$$
In this situation we take $\Pi(x,y)=x$.
The corresponding (Isaacs) equation is now
$$
\check{H}[\check{v}](x,y):=(1/2)\frac{\partial^{2}}{(\partial x)^{2}}
\check{v}(x,y)- y\frac{\partial^{2}}{ \partial x \partial y }
\check{v}(x,y)
+(1/2)y^{2}\frac{\partial^{2}}{(\partial y)^{2}}
\check{v}(x,y)
$$
$$
+\frac{\partial }{ \partial x}\check{v}(x,y)+yf(x)=0.
$$ 
As $G(x)$ and $\check{G}(x,y)$ one can take $1-|x|^2$ and
set $r(x,y)=y$.

It is a trivial computation to show that if $u(x)$
satisfies $H[u](x)\leq0$ at a point $x\in(-1,1)$, then for $\check{u}(x,y)
:=yu(x)$ we have $\check{H}[\check{u}](x,y)\leq0$   
for any $y>0$ and if we reverse the sign of the 
first inequality the same will happen with the second one.

By Theorem \ref{theorem 2.19.1} we have that
$\check{v}(x,y)=yv(x)$ in $D^{\v}_{\varepsilon}$ and since for $y=1$
$$
y_{t}=e^{-w_{t}-(1/2)t},
$$
we conclude that for any $\varepsilon\in(0,1)$
\begin{equation}
                                                \label{3.18.7}
v(x)=E\big[\int_{0}^{\tau^{\varepsilon}_{x} }e^{-w_{t}-(1/2)t}
f(x+w_{t}+t)\,dt
+y_{\tau^{\varepsilon}_{x }}v(x_{\tau^{\varepsilon}_{x} })\big],
\end{equation}
where $\tau^{\varepsilon}_{x}$ is the minimum
of the first exit time of $x+w_{t}+t$ from $(-1,1)$
and the first exit time of $e^{-w_{t}-(1/2)t}$ from
$(\varepsilon,\varepsilon^{-1})$. The latter
tends to infinity as $\varepsilon\downarrow0$
and we obtain \eqref{3.18.6} from \eqref{3.18.7}
and the fact that $v=0$ at $\pm1$.

The reader might have noticed that
the process given by \eqref{3.18.5}
is degenerate. It shows 
why in Assumption \ref{assumption 3.25.1}
we require only $\Pi(\check{x}_{t})$
 to be uniformly nondegenerate.

\end{remark}

\mysection{Main results for the whole space}
                                        \label{section 5.31.1}

In this section we keep the assumptions of Section
\ref{section 2.26.3} apart from Assumptions \ref{assumption 3.19.1}
and \ref{assumption 5.24.1} concerning the existence
of the barrier functions $G$ and $\check G$ and take
$D=\bR^{d}$. In case we encounter expressions
like $v(x_{\gamma})$ we set them to be zero
on the event $\{\gamma=\infty\}$.
In the whole space we need the following.
\begin{assumption}
                                         \label{assumption 5.31.1}

(i) The functions $ c,f, \check c,
\check f$ are bounded.

(ii) For a constant $\chi>0$ we have
$c^{\alpha\beta}(p,x ),\check
c^{\alpha\beta}(p,\check x )\geq\chi$ for all $\alpha,\beta,p,x$
and $\check x$.
\end{assumption}

Notice that in this situation $\tau^{\alpha_{\cdot}\beta_{\cdot}x}
=\infty$, however $\check \tau^{\alpha_{\cdot}\beta_{\cdot}
\check x}$ may still be finite.

\begin{theorem}
                                             \label{theorem 5.31.1}
Under the above assumptions all assertions of Theorems
\ref{theorem 1.14.1} and \ref{theorem 2.19.1} hold true.
\end{theorem}

Proof. First we deal with Theorem \ref{theorem 1.14.1}.
Take $D=D_{n}=\{x:|x|<n\}$ and $0$
in the original Theorem \ref{theorem 1.14.1}
in place of $D$ and $g$, respectively,
and  denote thus obtained function $v$ by $v_{n}$.
It is not hard to check that, due to the boundedness
of $f$ and the condition that $c\geq\chi$, in any
compact set $\Gamma\subset\bR^{d}$ we have
$v_{n}\to v$ uniformly on $\Gamma$ as $n\to\infty$.
Furthermore, since the boundary of $D_{n}$ is smooth
and $\sigma,b,c$ are bounded and $a$ is uniformly nondegenerate,
for each $n$ there exists a global barrier $G_{n}$
satisfying Assumption \ref{assumption 3.19.1} with $D_{n}$
in place of $D$. Therefore, 
by Theorem \ref{theorem 1.14.1}, $v_{n}$ are continuous and so is $v$.

For each $n\geq m\geq1$ we also have by Theorem \ref{theorem 1.14.1}
that
$$
v_{n}(x)=\infsup_{\bbeta\in\bB\,\,\alpha_{\cdot}\in\frA}
E_{x}^{\alpha_{\cdot}\bbeta(\alpha_{\cdot})}\big[
v_{n}(x_{\gamma\wedge\tau_{m}})e^{-\phi_{\gamma\wedge\tau_{m}}
-\psi_{\gamma\wedge\tau_{m}}}
$$
$$
+\int_{0}^{\gamma\wedge\tau_{m}}
\{f( x_{t})+\lambda_{t}v_{n}(x_{t})\}e^{-\phi_{t}-\psi_{t}}\,dt \big],
$$
where $\tau_{m}^{\alpha_{\cdot}\beta_{\cdot}x}$
is the first exit time of $x_{t}^{\alpha_{\cdot}\beta_{\cdot}x}$
from $D_{m}$. Since $v_{n}\to v$ uniformly on $\bar{D}_{m}$,
we conclude that
$$
v (x)=\infsup_{\bbeta\in\bB\,\,\alpha_{\cdot}\in\frA}
E_{x}^{\alpha_{\cdot}\bbeta(\alpha_{\cdot})}\big[
v (x_{\gamma\wedge\tau_{m}})e^{-\phi_{\gamma\wedge\tau_{m}}
-\psi_{\gamma\wedge\tau_{m}}}
$$
$$
+\int_{0}^{\gamma\wedge\tau_{m}}
\{f( x_{t})+\lambda_{t}v (x_{t})\}e^{-\phi_{t}-\psi_{t}}\,dt \big].
$$
Passing to the limit as $m\to\infty$ proves our theorem
in what concerns Theorem \ref{theorem 1.14.1}.

In case of Theorem \ref{theorem 2.19.1} the argument
is quite similar and we only comment on the existence
of $\check G_{n}$ satisfying Assumption
\ref{assumption 5.24.1} with $D^{\v}_{n}=
\{\check x\in D^{\v}:\Pi(\check x)\in D_{n}\}$.
Under obvious circumstances one can take
 $\check G_{n}(\check x) =G_{n}(\Pi(\check x))$.
In the general case one should construct $G_{n}$
for operators with, perhaps, a smaller ellipticity constant
and larger drift terms.
The theorem is proved.

\mysection{An auxiliary result}
                                               \label{section 2.27.1}

In this section $D$ is not assumed to be bounded.
We need a bounded continuous function $\Psi$ on $\bar{D}$ such that
 $\Psi\geq0$ in $D$
and $\Psi=0$ on $\partial D$ (if $\partial D\ne\emptyset$).
We assume that we are given two continuous
$\cF_{t}$-adapted processes $x'_{t}$ and $x''_{t}$ in $\bR^{d}$ with
$x'_{0},x''_{0}\in D$ (a.s.) and progressively measurable 
real-valued processes
$c'_{t},c''_{t},f'_{t},f''_{t}$.  Suppose that $c',c''\geq0$.

Define $\tau'$ and $\tau''$ as the
first exit times of $x'_{t}$ and $x''_{t}$ from $D$,
respectively. Then  introduce
$$
\phi'_{t}=\int_{0}^{t}c'_{s}\,ds,\quad
\phi''_{t}=\int_{0}^{t}c''_{s}\,ds,
$$
and suppose that
\begin{equation}
                                              \label{6.14.1}
E \int_{0}^{ \tau'}|f'_{t}|e^{-\phi'_{t}}\,dt+
E \int_{0}^{\tau'' }|f''_{t}|e^{-\phi''_{t}}\,dt<\infty.
\end{equation}
 
\begin{remark}
                                           \label{remark 2.26.3}
According to Theorem 2.2.1  of \cite{Kr77} the above 
requirements about $f$ and $c$ 
are fulfilled if Assumption \ref{assumption 1.9.1} is satisfied
and we take $x_{t}$ and $(f, c)$ with prime and double prime
of the type
$$
x^{\alpha_{\cdot}\beta_{\cdot}x}_{t},\quad
(f,c)^{\alpha_{t}\beta_{t}}(
x^{\alpha_{\cdot}\beta_{\cdot}x}_{t}
),
$$
respectively, where $\alpha_{\cdot}\in\frA$,
 $\beta_{\cdot}\in\frB$, and $x\in\bR^{d}$.

\end{remark}
Finally set
$$
v'=E \int_{0}^{\tau'}f'_{t}e^{-\phi'_{t}}\,dt,
\quad
v''=E \int_{0}^{\tau''}f''_{t}e^{-\phi''_{t}}\,dt.
$$
 
 Now comes our main assumption.
\begin{assumption}
                                        \label{assumption 2.24.1}
The processes 
$$
\Psi(x'_{t\wedge\tau'})e^{-\phi'_{t\wedge\tau'}}
+\int_{0}^{t\wedge\tau'}e^{-\phi'_{s}}\,ds ,\quad
\Psi(x''_{t\wedge\tau''})e^{-\phi''_{t\wedge\tau''}}
+\int_{0}^{t\wedge\tau''}e^{-\phi''_{s}}\,ds 
$$
are supermartingale.

\end{assumption}
\begin{remark}
                                              \label{remark 2.26.2}
Observe that Assumption \ref{assumption 2.24.1}
is satisfied under the assumptions of Theorem  
\ref{theorem 1.14.1} if we take $\Psi=G$ from Theorem \ref{theorem 1.14.1}
and other objects from Remark \ref{remark 2.26.3}.

Indeed, by It\^o's formula
$$
 G(x _{t\wedge\tau})e^{-\phi_{t\wedge\tau}}
+\int_{0}^{t\wedge\tau}e^{-\phi_{s}}\,ds
$$
is a local supermartingale,  
where 
\begin{equation}
                                                          \label{4.6.7}
x_{t}=x^{\alpha_{\cdot}\beta_{\cdot}x}_{t},\quad
\tau=\tau^{\alpha_{\cdot}\beta_{\cdot}x},\quad
\phi_{t}=\phi^{\alpha_{\cdot}\beta_{\cdot}x}_{t}  .
\end{equation}
Since it is nonnegative or constant, it is a supermartingale.
\end{remark}

Denote
$$
\Phi_{t}=e^{-\phi'_{t}}+e^{-\phi''_{t}},\quad
\Delta_{c}=E\int_{0}^{\tau'\wedge\tau''}
|c'_{t} -
c''_{t} |\Phi_{t}\,dt,
$$
and by replacing $c$ with $f$ define $\Delta_{f}$.

\begin{lemma}
                                         \label{lemma 2.24.1}

Introduce a constant $M_{f}$ (perhaps
$M_{f}=\infty$) such that for each $t\geq0$ (a.s.)
\begin{equation}
                                                      \label{2.22.5}
I_{\tau'\wedge\tau''>t}E\big\{\int_{t}^{\tau'\wedge\tau''}
|f''_{s}| \Phi_{s}\,ds\mid
\cF_{t}\big\}\leq \Phi_{t} M_{f}.
\end{equation}
  Then
$$
|v'-
v''|
\leq \Delta_{f}+M_{f}\Delta_{c}
+ \sup |f'| EI_{\tau''<\tau'}[\Psi(x'_{ \tau''})-\Psi(x'' _{ \tau''})]
e^{-\phi'_{ \tau''}}
$$
\begin{equation}
                                                     \label{2.22.3}
+ \sup |f''| EI_{\tau'<\tau''}[\Psi(x''_{ \tau'})-\Psi(x' _{ \tau '})]
e^{-\phi''_{ \tau'}} ,
\end{equation}
where   the last two terms
can be dropped if $\tau'=\tau''$ (a.s.).

\end{lemma}

Proof. We have  
$$
\big|v''
-E\int_{0}^{\tau'\wedge \tau''}
f''_{t} e^{-\phi''_{t}}\,dt\big|\leq
\sup|f''|E\int_{\tau'\wedge\tau''}^{\tau''}e^{-\phi''_{t}}\,dt,
$$
where owing to \eqref{6.14.1},
Assumption \ref{assumption 2.24.1}, and the fact that bounded
$\Psi\geq0$,
 the last expectation is
dominated by
$$
E\Psi(x''_{\tau'\wedge\tau''})e^{-\phi''_{\tau'\wedge\tau''}}
I_{\tau'<\tau'' }=
EI_{\tau'<\tau''}\Psi(x''_{\tau' })e^{-\phi''_{ \tau '}}
$$
$$
=EI_{\tau'<\tau''}[\Psi(x''_{ \tau'})-\Psi(x' _{ \tau '})]
e^{-\phi''_{ \tau'}}.
$$
Similar estimates hold for $v'$ and this
 shows how the last terms in \eqref{2.22.3}
appear and when they disappear.

Next,
$$
E\int_{0}^{\tau'\wedge\tau''}\big|
f'_{t} e^{-\phi'_{t}} 
-
f''_{t} e^{-\phi''_{t}}\big|\,dt
\leq \Delta_{f}+J,
$$
where
$$
J=E\int_{0}^{\tau'\wedge\tau''}
|f''_{t}|\,| e^{-\phi''_{t}}
-e^{-\phi' _{t}}|\,dt
\leq E\int_{0}^{\tau'\wedge\tau''}
|f''_{t}|C_{t}
\Phi_{t}\,dt,
$$
$$
C_{t}=\int_{0}^{t}|c'_{s}-
c''_{s}|\,ds.
$$
By using Fubini's theorem it is easily seen that
the last expectation above equals
$$
E\int_{0}^{\tau'\wedge\tau''}\big(\int_{s}^{\tau'\wedge\tau''}
|f''_{t}| 
\Phi_{t}\,dt\big)
|c'_{s}-
c''_{s}|\,ds,
$$
which owing to \eqref{2.22.5} is less than $M_{f}\Delta_{c}$.
This proves the lemma.

\begin{remark}
                                         \label{remark 2.22.3}
Assumption \eqref{2.22.5} is satisfied if, for instance,
for each $t\geq 0$
\begin{equation}
                                                    \label{4.16.6}
I_{\tau''>t}E\big\{\int_{t}^{\tau''}
|f''_{s}|  \,ds\mid
\cF_{t}\big\}\leq   M_{f}.
\end{equation}
Indeed, in that case the left-hand side of 
\eqref{2.22.5} is less that $\Phi_{t}$ times the left-hand side of
\eqref{4.16.6} just because $\Phi_{t}$ is a decreasing function
of $t$.

This observation will be later used in conjunction with
Theorem 2.2.1  of \cite{Kr77}. 
\end{remark}

\mysection{A general approximation result from above}
                                              \label{section 2.25.1} 

In this section Assumption \ref{assumption 1.9.1} (iv) about the
uniform nondegeneracy as well as Assumption
\ref{assumption 3.19.1} concerning  $G$ 
are not used  and the domain $D$ 
  is not supposed to be bounded.

We impose the following.
\begin{assumption}
                                           \label{assumption 3.13.1}
(i) Assumptions \ref{assumption 1.9.1} (i) b), (ii)  
are satisfied.

(ii) The functions $c^{\alpha\beta}( x)$ and $f^{\alpha\beta}( x)$
are bounded on $A\times B \times \bR^{d}$ and uniformly
 continuous with respect to $x\in\bR^{d}$
uniformly with respect to $\alpha,\beta $.

\end{assumption}  

Set 
$$
A_{1}=A 
$$
 and let $A_{2}$ be a 
separable metric space having no common points with $A_{1}$.

\begin{assumption}
                                         \label{assumption 4.29.2}
The functions $ 
\sigma^{\alpha\beta}( x)$,
$ b^{\alpha\beta}( x)$, $ 
c^{\alpha\beta}( x)$, and
$ f^{\alpha\beta}( x)$ 
 are also defined on
  $A_{2}\times B \times\bR^{d}$ in such a way that they are
{\em independent\/} of $\beta$ and 
  Assumptions \ref{assumption 1.9.1} (i) b), (ii) are satisfied
with, perhaps,   larger constants $K_{0}$  $K_{1}$
and, of course, with $A_{2}$ in place of $A$.
The functions $ 
c^{\alpha\beta}( x)$  and
$ f^{\alpha\beta}( x)$ are bounded on $A_{2}\times B
 \times \bR^{d}$.
\end{assumption}

Define
$$
\hat{A}=A_{1}\cup A_{2}.
$$

Then we introduce $\hat{\frA}$ as the set of progressively measurable
$\hat{A}$-valued processes and $\hat{\bB}$ as the set
of $\frB$-valued functions $ \bbeta(\alpha_{\cdot})$
on $\hat{\frA}$ such that,
for any $T\in[0,\infty)$ and any $\alpha^{1}_{\cdot},
\alpha^{2}_{\cdot}\in\hat{\frA}$ satisfying
$$
P(  \alpha^{1}_{t}=\alpha^{2}_{t} 
 \quad\text{for almost all}\quad t\leq T)=1,
$$
we have
$$
P(  
 \bbeta_{t}(\alpha^{1}_{\cdot})=\bbeta_{t}(\alpha^{2}_{\cdot}) 
 \quad\text{for almost all}\quad t\leq T)=1.
$$

\begin{assumption}
                                           \label{assumption 2.23.1}
There exists a bounded uniformly continuous
in $\bar{D}$ function $G\in C^{2}_{loc}(D)$ such that
$G=0$ on $\partial D$ (if $D\ne\bR^{d}$) and
$$
L^{\alpha\beta}G( x)\leq-1
$$
in $  D$ for all $\alpha\in\hat{A}$ and $\beta\in B$.
\end{assumption}

Here are a few consequences of Assumption \ref{assumption 2.23.1}.

\begin{lemma}
                                         \label{lemma 2.23.1}
For any constant $\chi\leq (2\sup_{D}G)^{-1}$ and any
$\alpha_{\cdot}\in\hat{\frA}$, $\beta_{\cdot}\in\frB$,
and $x\in\bar{D}$ the process
$$
 G(x _{t\wedge
\tau})
e^{ \chi(t\wedge\tau )
-\phi _{t\wedge\tau }}
+(1/2)\int_{0}^{t\wedge\tau }
e^{ \chi s-\phi _{s}}\,ds,
$$
where we use notation \eqref{4.6.7},
is a supermartingale and
$$
E^{\alpha_{\cdot}\beta_{\cdot}}_{x}\int_{0}^{\tau}
e^{\chi t-\phi_{t}}\,dt\leq 2G(x).
$$
In particular, for any $T\in[0,\infty)$
$$
E^{\alpha_{\cdot}\beta_{\cdot}}_{x}I_{\tau>T}\int_{T}^{\tau}e^{ -\phi_{t}}\,dt
= e^{-\chi T}E^{\alpha_{\cdot}\beta_{\cdot}}_{x}I_{\tau>T}\int_{T}^{\tau}e^{\chi T -\phi_{t}}\,dt
$$
$$
\leq e^{-\chi T}E^{\alpha_{\cdot}\beta_{\cdot}}_{x}\int_{0}^{\tau}
e^{\chi t-\phi_{t}}\,dt
\leq 2e^{-\chi T}G(x).
$$

Finally, for any stopping time $\gamma\leq
\tau^{\alpha_{\cdot}\beta_{\cdot} x}$
$$
E^{\alpha_{\cdot}\beta_{\cdot}}_{x}I_{\gamma>T}
G(x_{\gamma})e^{-\phi_{\gamma}}\leq
E^{\alpha_{\cdot}\beta_{\cdot}}_{x}I_{\gamma>T}
G(x_{T})e^{-\phi_{T}}
$$
$$
\leq e^{-\chi T}
E^{\alpha_{\cdot}\beta_{\cdot}}_{x} I_{\gamma>T}
G(x_{T})e^{\chi T-\phi_{T}}
$$
$$
\leq e^{-\chi T}
E^{\alpha_{\cdot}\beta_{\cdot}}_{x}  
G(x_{T\wedge\gamma})e^{\chi(T\wedge\gamma)-
\phi_{T\wedge\gamma}}\leq e^{-\chi T}G(x).
$$
 
\end{lemma}

The proof of this lemma is easily achieved by using
It\^o's formula and the fact that $L^{\alpha\beta}G+\chi G
\leq-1/2$ on $D$ for all $\alpha,\beta$.

Take a constant $K\geq0$ and set 
$$
v_{K}(x)=\infsup_{\bbeta\in\hat{\bB}\,\,\alpha_{\cdot}\in\hat{\frA}}
v^{\alpha_{\cdot}\bbeta(\alpha_{\cdot})}_{K}(x),
$$
where
$$
v^{\alpha_{\cdot}\beta _{\cdot} }_{K}(x)=
E_{x}^{\alpha_{\cdot}\beta _{\cdot} }\big[\int_{0}^{\tau}
 f_{K} ( x_{t})e^{-\phi_{t} }\,dt
+g(x_{\tau})e^{-\phi_{\tau}
 }\big]
$$
$$
= :v^{\alpha_{\cdot}\beta _{\cdot} }  (x)-K
E_{x}^{\alpha_{\cdot}\beta _{\cdot} }\int_{0}^{\tau}
I_{\alpha_{t}\in A_{2}}e^{-\phi_{t} }\,dt,
$$
$$
 f^{\alpha\beta}_{K}( x)=f^{\alpha\beta}( x)-KI_{\alpha\in A_{2}}.
$$
Observe that
$$
v(x)=\infsup_{\bbeta\in\bB\,\,\alpha_{\cdot}\in\frA}
v^{\alpha_{\cdot}\bbeta(\alpha_{\cdot})}(x).
$$

These definitions make sense  
owing to Lemma \ref{lemma 2.23.1}, which also implies that
  $v^{\alpha_{\cdot}\beta _{\cdot} }_{K}$
and $v^{\alpha_{\cdot}\beta _{\cdot} }$  
and bounded in $\bar{D}$.

\begin{theorem}
                                             \label{theorem 1.20.1}
We have
 $v_{K}\to v$ uniformly on $\bar{D}$ as $K\to\infty$.
\end{theorem}

We need the following.

\begin{lemma}
                                           \label{lemma 3.24.1}
There exists a constant $N$ depending
only on $K_{0},K_{1}$,  and $d$ 
such that for any
$\alpha_{\cdot}\in\hat{\frA}$, $\beta_{\cdot}\in\frB$,
 $x\in\bR^{d}$,  $T\in[0,\infty)$, and
stopping time $\gamma$
$$
E^{\alpha_{\cdot}\beta_{\cdot}}_{x}
\sup_{t\leq T\wedge\gamma}|x_{t}-y_{t}|\leq
Ne^{NT} \big(E^{\alpha_{\cdot}\beta_{\cdot}}_{x}\int_{0}^{
T\wedge\gamma}I_{\alpha_{t}\in A_{2}}\,dt\big)^{1/2} ,
$$
where 
$$
y^{\alpha_{\cdot} \beta_{\cdot}x}_{t}
=x^{\pi\alpha_{\cdot} \beta_{\cdot}x}_{t}.
$$
\end{lemma}

Proof. For simplicity of notation we drop the
superscripts $\alpha_{\cdot} ,\beta_{\cdot},x$.
Observe that $x_{t}$ and $y_{t}$ satisfy
$$
x_{t}=x+\int_{0}^{t}\sigma^{\alpha_{s}\beta_{s}}( x_{s})
\,dw_{s}+\int_{0}^{t}b^{\alpha_{s}\beta_{s}}( x_{s})
\,ds,
$$
$$
y_{t}=x+\int_{0}^{t}\sigma^{\alpha_{s}\beta_{s}}( y_{s})
\,dw_{s}+\int_{0}^{t}b^{\alpha_{s}\beta_{s}}( y_{s})
\,ds+\eta_{t},
$$
where $\eta_{t}=I_{t}+J_{t}$,
$$
I_{t}=\int_{0}^{t}[\sigma^{\pi\alpha_{s}\beta_{s}}( y_{s})
-\sigma^{\alpha_{s}\beta_{s}}( y_{s})]\,dw_{s},
$$
$$
J_{t}=\int_{0}^{t}[b^{\pi\alpha_{s}\beta_{s}}( y_{s})
-b^{\alpha_{s}\beta_{s}}( y_{s})]\,ds.
$$

By Theorem II.5.9 of \cite{Kr77} (where we replace
the processes $x_{t}$ and $\tilde{x}_{t}$ with appropriately stopped ones)
for any $T\in[0,\infty)$ and any stopping time $\gamma
 $
\begin{equation}
                                                       \label{1.22.04}
E\sup_{t\leq T\wedge \gamma}|x _{t}-y _{t}|^{2}\leq
Ne^{NT}E\sup_{t\leq T\wedge \gamma}|I _{t}+J _{t}|^{2},
\end{equation}
where $N$ depends only on $K_{1}$  and $d$, 
which by Theorem III.6.8
of \cite{Kr95} leads to
\begin{equation}
                                                       \label{1.23.01}
E\sup_{t\leq T\wedge\gamma }|x _{t}-y _{t}|
\leq Ne^{NT}E\sup_{t\leq T\wedge\gamma }|I _{t}+J _{t}|
\end{equation}
with the constant $N$ being three times the one from \eqref{1.22.04}.

By using Davis's inequality   we see that for any $T\in[0,\infty)$
$$
E\sup_{t\leq T\wedge\gamma }|I _{t}| 
\leq NE\big(\int_{0}^{T\wedge\gamma}
I_{\alpha _{s}\in A_{2}}\,ds\big)^{1/2}
\leq N\big(E\int_{0}^{T\wedge\gamma}
I_{\alpha _{s}\in A_{2}}\,ds\big)^{1/2}.
$$

Furthermore, almost obviously
$$
E\sup_{t\leq T\wedge\gamma }|J_{t}| 
\leq N E\int_{0}^{T\wedge\gamma}
I_{\alpha _{s}\in A_{2}}\,ds 
\leq NT^{1/2}\big(E\int_{0}^{T\wedge\gamma}
I_{\alpha _{s}\in A_{2}}\,ds\big)^{1/2}
$$
and this in combination with \eqref{1.23.01} proves the lemma.

{\bf Proof of Theorem \ref{theorem 1.20.1}}.
Without losing generality we may
assume  that $g\in C^{3}(\bR^{d})$ since
the functions of this class   uniformly  
approximate any $g$ which is uniformly continuous 
in $\bR^{d}$. Then notice that by It\^o's formula
and Lemma \ref{lemma 2.23.1}
for $g\in C^{3}(\bR^{d})$  
we have
$$
E_{x}^{\alpha_{\cdot}\beta _{\cdot} }\big[\int_{0}^{\tau}
 f_{K} ( x_{t})e^{-\phi_{t} }\,dt
+g(x_{\tau})e^{-\phi_{\tau}
 }\big]
$$
$$
=g(x)+E_{x }^{\alpha _{\cdot} \beta 
 _{\cdot} } \int_{0}^{\tau}[ 
 \hat{f} ( x_{t} )-KI_{\alpha_{t}\in A_{2}}] e^{-\phi_{t} }\,dt,
$$
where 
$$
\hat{f} ^{\alpha\beta}( x ):=f^{\alpha\beta} ( x ) 
+ L^{\alpha\beta}g( x),
$$
which
is bounded and, for $(\alpha,\beta)\in A\times B$, is
 uniformly continuous in $x$ uniformly with respect to $\alpha,\beta$.
This argument   shows that without losing generality
we may (and will) also assume that $g=0$.

Next, since   $\frA\subset\hat{\frA}$ and for $\alpha_{\cdot}\in
\hat{\frA}$
and $\bbeta \in\hat{\bB}$ we have $\bbeta(\alpha_{\cdot})\in\frB$,
it holds that
 $$
v_{K}\geq v.
$$

To estimate $v_{K}$ from above,
we need a mapping $\pi:\hat{A}\to A_{1}$ defined as 
$\pi(\alpha)=\alpha$ if $\alpha\in A_{1}$ and
$\pi(\alpha)=\alpha^{*}$ if $\alpha\in A_{2}$, where
$\alpha^{*}$ is a fixed point in $A$.
Take $\bbeta\in\bB$ and define $\hat{\bbeta} \in\hat{\bB}$ by
$$
\hat{\bbeta}_{t}(\alpha_{\cdot})=\bbeta_{t}(\pi\alpha_{\cdot}).
$$
Also take any sequence $x^{n}
 \in \bar{D}$, $n=1,2,...$, and find a sequence
$\alpha^{n}_{\cdot}\in\hat{\frA}$ such that
$$
v_{K}(x^{n}) \leq  \sup_{\alpha\in\hat{\frA}}
E_{x^{n}}^{\alpha_{\cdot}\hat{\bbeta}(\alpha_{\cdot})}
 \int_{0}^{\tau}
 f _{K} ( x_{t})e^{-\phi_{t} }\,dt
$$
\begin{equation}
                                                    \label{2.22.6}
=1/n+v^{\alpha^{n}_{\cdot}\hat{\bbeta}(\alpha^{n}_{\cdot})}(x^{n})
-KE\int_{0}^{\tau^{n} }I_{\alpha^{n}_{t}\in A_{2}}
e^{-\phi^{n}_{t} }\,dt,
\end{equation}
where
$$
(\tau^{n},\phi^{n}_{t} )
=( \tau,\phi _{t} )
^{\alpha^{n}_{\cdot}\hat{\bbeta}(\alpha^{n}_{\cdot})x^{n}}.
$$

It follows from  Lemma \ref{lemma 2.23.1}   that
 there is a constant $N$
independent of $n$ and $K$ such that $|v^{\alpha^{n}_{\cdot}\hat{\bbeta}
(\alpha^{n}_{\cdot})}(x^{n})|\leq N$, $|V|\leq N$, $v_{K}\geq v\geq-N$
and we conclude from \eqref{2.22.6} that for any $T\in[0,\infty)$
and 
$$
\bar{c}:=\sup c 
$$
 we have
\begin{equation}
                                                    \label{1.20.2}
 E\int_{0}^{\tau^{n} }I_{\alpha^{n}_{t}\in A_{2}}
e^{-t\bar{c}}\,dt\leq N/K,\quad
 E\int_{0}^{\tau^{n} \wedge T}I_{\alpha^{n}_{t}\in A_{2}}
 \,dt\leq Ne^{NT}/K,
\end{equation}
where and below in the proof by $N$ we denote constants
which may change from
one occurrence to another and independent of $n$, $K$, and $T$.

Next, introduce
$$
x^{n}_{t} =
x^{\alpha^{n}_{\cdot}\hat{\bbeta} 
(\alpha^{n}_{\cdot}) x^{n}}_{t} ,\quad
 y^{n}_{t} 
=x^{\pi\alpha^{n}_{\cdot}\hat{\bbeta} 
( \alpha^{n}_{\cdot}) x^{n} }_{t} ,\quad
\pi\phi^{n}_{t}=\int_{0}^{t}c^{ \pi\alpha^{n}_{s}\hat{\bbeta} _{s}
( \alpha^{n}_{\cdot})}(y^{n}_{s} )\,ds,
$$
define $\gamma^{n}$ as the first exit time of $y^{n}_{t}$
from $D$,
and, with the aim of applying Lemma \ref{lemma 2.24.1},
observe that by identifying $x^{n}_{t},y^{n}_{t},\tau^{n},\gamma^{n}$
 and the objects
related to them with $x'_{t},x''_{t},\tau',\tau''$ and the objects related
to them, respectively,
  we have
$$
|c'_{t}-c''_{t}|
=|c^{\alpha^{n}_{t}\hat{\bbeta}_{t}(\alpha^{n}_{\cdot})}
( x^{n}_{t})-
c^{\pi\alpha^{n}_{t}\hat{\bbeta}_{t}(\alpha^{n}_{\cdot})}
( y^{n}_{t})| .
$$
Hence
 for any $T\in(0,\infty)$  
$$
\Delta_{c}^{n} =E\int_{0}^{\tau^{n}\wedge\gamma^{n}}
|c^{\alpha^{n}_{t}\hat{\bbeta}_{t}(\alpha^{n}_{\cdot})}
( x^{n}_{t})-c^{\pi\alpha^{n}_{t}\hat{\bbeta}_{t}(\alpha^{n}_{\cdot})}
( y^{n}_{t})|(e^{-\phi_{t}^{n} }+e^{-\pi\phi_{t}^{n}
 })\,dt
$$
$$
\leq E\int_{0}^{\tau^{n}\wedge\gamma^{n}\wedge T} W_{c }
(|x^{n}_{t}
-y^{n}_{t}|) \,dt+I_{n }+J_{n },
$$
where $W_{c }$ is the modulus of continuity of $c $
 and
$$
I_{n}=N  E\int_{0}^{\tau^{n}\wedge\gamma^{n}\wedge T}
I_{\alpha^{n}_{t}\in A_{2}}
 \,dt,
$$
$$
J_{n }=N E\int_{\tau^{n}\wedge\gamma^{n}\wedge T}^{\tau^{n}\wedge\gamma^{n}}
(e^{-\phi_{t}^{n}}+e^{-\pi\phi_{t}^{n}})\,dt.
$$
By virtue of \eqref{1.20.2} we have $I_{n}\leq Ne^{NT }/K$
and 
 $
J_{n }\leq N e^{-\chi T} 
 $
by Lemma \ref{lemma 2.23.1}, say with $\chi= (2\sup_{D}G)^{-1}$.  
Therefore,
$$
\Delta_{c}^{n}\leq TEW_{c}(\sup_{t\leq\tau^{n}\wedge T}
|x^{n}_{t}-y^{n}_{t}|)+ Ne^{NT}/K+Ne^{-\chi T}.
$$
A similar estimate holds if we replace $c$ with $f$.

As long as the last terms in \eqref{2.22.3} are concerned, observe  
that  
$$
 E|G(  x^{n}_{\tau^{n} })-G(
y^{n}_{\tau^{n} } )|
 e^{-\pi\phi^{n}_{\tau^{n}}}I_{\tau^{n}<\gamma^{n}}
$$
$$
\leq EW_{G}\big(\sup_{t\leq\tau^{n}\wedge\gamma^{n} \wedge T}
|x^{n}_{t}-y^{n}_{t}|\big)+R_{n},
$$
where $W_{G}$ is the modulus of continuity
of $G$ and 
$$
R_{n}=EI_{\gamma^{n}>\tau^{n}>T}G(y^{n}_{\tau^{n}})
e^{-\pi\phi_{\tau^{n}}}
\leq EI_{\gamma^{n}\wedge\tau^{n}>T}G(y^{n}_{\gamma^{n}\wedge\tau^{n}})
e^{-\pi\phi_{\gamma^{n}\wedge\tau^{n}}} \leq Ne^{-\chi T},
$$
with the second inequality following from Lemma \ref{lemma 2.23.1}.
 
Finally, in light of Lemma \ref{lemma 2.23.1} one can take
$M_{f}$ in  Lemma \ref{lemma 2.24.1}
to be a constant $N$ independent of $n$ and $K$
and then by applying
  Lemma \ref{lemma 2.24.1} we conclude from
\eqref{2.22.6} that
$$
v_{K}(x^{n})\leq 1/n
+v^{\pi\alpha_{\cdot}^{n}\bbeta(\pi\alpha_{\cdot}^{n})}(x^{n})
$$
$$
+(T+1)EW (\sup_{t\leq\tau^{n}\wedge\gamma^{n}\wedge T}
|x^{n}_{t}-y^{n}_{t}|)+ Ne^{NT }/K+Ne^{-\chi T},
$$
where $W(r)$ is a bounded function such that $W(r)\to0$ as $r\downarrow
0$.

  This result, \eqref{1.20.2},
and Lemma \ref{lemma 3.24.1} imply that,
for any $T$,
\begin{equation}
                                                      \label{1.25.1}
v_{K}(x^{n})\leq 1/n
+v^{\pi\alpha_{\cdot}^{n}\bbeta(\pi\alpha_{\cdot}^{n})}(x^{n})
+w(T,K)+ Ne^{NT }/K+Ne^{-\chi T},
\end{equation}
where 
 $w(T,K)$ is independent of   $n$ and
 $w(T,K)\to0$ as $K\to\infty$ for any $T$.
Hence
$$
v_{K}(x^{n})\leq \sup_{\alpha_{\cdot}\in\frA}
 v^{ \alpha _{\cdot}
 \bbeta ( \alpha_{\cdot} )}(x^{n}) 
+w(T,K)+ Ne^{NT }/K+Ne^{-\chi T}+1/n.
$$
Owing to  the arbitrariness of $\bbeta\in\bB$ we have
$$
v_{K}(x^{n})\leq v(x^{n})+w(T,K)+ Ne^{NT }/K+Ne^{-\chi T}+1/n,
$$
and  the arbitrariness of $x^{n}$ yields
$$
\sup_{\bar{D}}(v_{K}-v) \leq  w(T,K)+ Ne^{NT }/K+Ne^{-\chi T} ,
$$
which leads to the desired result after first letting $K\to\infty$
and then $T\to\infty$. The theorem is proved.

\mysection{A general approximation result from below}
                                              \label{section 3.18.1}

 As in Section \ref{section 2.25.1},  
  Assumption \ref{assumption 1.9.1} (iv) about the  
uniform nondegeneracy as well as Assumption
\ref{assumption 3.19.1} concerning  $G$ 
are not used  and the domain $D$ 
  is not supposed to be bounded.

However,
we  suppose that Assumption \ref{assumption 3.13.1}
is satisfied.
Here we allow $\beta$ to change in a larger set
penalizing using controls other than initially
available. 

Set 
$$
B_{1}=B 
$$
 and let $B_{2}$ be a 
separable metric space having no common points with $B_{1}$.
\begin{assumption}
                                         \label{assumption 4.29.1}
The functions $ 
\sigma^{\alpha\beta}( x)$,
$ b^{\alpha\beta}(x)$, $ 
c^{\alpha\beta}(x)$, and
$ f^{\alpha\beta}(x)$ 
 are also defined on  $A\times B_{2} \times\bR^{d}$ in 
such a way that they are
{\em independent\/} of $\alpha$   and 
 Assumptions \ref{assumption 1.9.1} (i) b), (ii) are satisfied
with, perhaps,  larger constants $K_{0}$  and $K_{1}$
and, of course, with $B_{2}$ in place of $B$.
The functions $ 
c^{\alpha\beta}(x)$  and
$ f^{\alpha\beta}(x)$ are bounded on $A\times B_{2}
 \times\bR^{d}$.
\end{assumption}

Define
$$  
\hat{B}=B_{1}\cup B_{2}. 
$$

Then we introduce $\hat{\frB}$ as the set of progressively measurable
$\hat{B}$-valued processes and $\hat{\bB}$ as the set
of $\hat{\frB}$-valued functions $ \bbeta(\alpha_{\cdot})$
on $ \frA $ such that,
for any $T\in[0,\infty)$ and any $\alpha^{1}_{\cdot},
\alpha^{2}_{\cdot}\in \frA $ satisfying
$$
P( \alpha^{1}_{t}=\alpha^{2}_{t} \quad\text{for almost all}\quad t\leq T)=1,
$$
we have
$$
P(  \bbeta_{t}(\alpha^{1}_{\cdot})=\bbeta_{t}(\alpha^{2}_{\cdot})
\quad\text{for almost all}\quad t\leq T)=1.
$$

\begin{assumption}
                                           \label{assumption 2.23.01}
There exists a bounded uniformly continuous
in $\bar{D}$ function $G\in C^{2}_{loc}(D)$ such that
$G=0$ on $\partial D$ (if $D\ne\bR^{d}$) and
$$
L^{\alpha\beta}G(x)\leq-1
$$
in $  D$ for all $\alpha\in A$ and $\beta\in \hat{B}$.
\end{assumption}

Take a constant $K\geq0$ and set 
$$
v_{-K}(x)=\inf_{\bbeta\in\hat{\bB}}\sup_{\alpha_{\cdot}\in \frA} 
v^{\alpha_{\cdot}\bbeta(\alpha_{\cdot})}_{-K}(x),
$$
where
$$
v^{\alpha_{\cdot}\beta _{\cdot} }_{-K}(x)=
E_{x}^{\alpha_{\cdot}\beta _{\cdot} }\big[\int_{0}^{\gamma\wedge\tau}
 f_{-K} ( x_{t})e^{-\phi_{t} }\,dt
+g(x_{\gamma\wedge\tau})e^{-\phi_{\gamma\wedge\tau}
 }\big]
$$
$$
= :v^{\alpha_{\cdot}\beta _{\cdot} }  (x)+K
E_{x}^{\alpha_{\cdot}\beta _{\cdot} }\int_{0}^{\gamma}
I_{\beta_{t}\in B_{2}}e^{-\phi_{t} }\,dt,
$$
$$
 f^{\alpha\beta}_{-K}( x)=f^{\alpha\beta}( x)
+KI_{\beta\in B_{2}} .
$$
We reiterate that
$$
v(x)=\infsup_{\bbeta\in\bB\,\,\alpha_{\cdot}\in\frA}
v^{\alpha_{\cdot}\bbeta(\alpha_{\cdot})}(x).
$$

These definitions make sense by the same reason as in Section
\ref{section 2.25.1}.  

\begin{theorem}
                                             \label{theorem 4.18.1}
We have
 $v_{-K}\to v$ uniformly on $\bar{D}$ as $K\to\infty$.
\end{theorem}

Proof. As in the proof of Theorem \ref{theorem 1.20.1}
we may assume that $g=0$. Then since $\bB\subset\hat{\bB}$ we have that
$v_{-K}\leq v$. To estimate $v_{-K}$ from below take any  sequence
$x^{n}\in\bar{D}$ and find a sequence $\bbeta^{n}\in\hat{\bB}$
such that
$$
v_{-K}(x^{n})\geq-1/n +\sup_{\alpha_{\cdot}\in\frA}
E^{\alpha_{\cdot}\bbeta^{n}(\alpha_{\cdot})}_{x^{n}}
\int_{0}^{\tau}f_{-K}( x_{t})e^{-\phi_{t} }
\,dt.
$$
Since the last supremum is certainly greater than a negative constant
independent of $n$ plus
$$
K\sup_{\alpha_{\cdot}\in\frA}
E^{\alpha_{\cdot}\bbeta^{n}(\alpha_{\cdot})}_{x^{n}}
\int_{0}^{\tau}I_{\bbeta_{t}(\alpha_{\cdot})
\in B_{2}}e^{-\bar{c}t}
\,dt,
$$
where $\bar{c}$ is the same as in Section \ref{section 2.25.1},
we conclude that
\begin{equation}
                                                    \label{3.18.3}
\sup_{\alpha_{\cdot}\in\frA}
E^{\alpha_{\cdot}\bbeta^{n}(\alpha_{\cdot})}_{x^{n}}
\int_{0}^{\tau}I_{\bbeta_{t}(\alpha_{\cdot})
\in B_{2}}e^{-\bar{c}t}
\,dt\leq N/K.
\end{equation}

Next, find a sequence of $\alpha^{n}_{\cdot}\in\frA$
such that
$$
E^{\alpha^{n}_{\cdot}\pi\bbeta^{n}(\alpha^{n}_{\cdot})}_{x^{n}}
\int_{0}^{\tau}f ( x_{t})e^{-\phi_{t} }
\,dt\geq v(x^{n})-1/n.
$$
By using \eqref{3.18.3} and arguing as in the proof 
of Theorem \ref{theorem 1.20.1}
one proves that
$$
I_{n}:=\big|E^{\alpha^{n}_{\cdot}\pi\bbeta^{n}(\alpha^{n}_{\cdot})}_{x^{n}}
\int_{0}^{ \tau}f ( x_{t})e^{-\phi_{t} }
\,dt
-
E^{\alpha^{n}_{\cdot} \bbeta^{n}(\alpha^{n}_{\cdot})}_{x^{n}}
\int_{0}^{ \tau}f ( x_{t})e^{-\phi_{t} }
\,dt\big| 
$$
tends to zero as $n\to\infty$. This leads to the desired result
since
$$
v_{-K}(x^{n})\geq  -1/n
+E^{\alpha^{n}_{\cdot}\bbeta^{n}(\alpha^{n}_{\cdot})}_{x^{n}}
\int_{0}^{ \tau}f ( x_{t})e^{-\phi_{t} }
\,dt
$$
$$
\geq -1/n+I_{n}+E^{\alpha^{n}_{\cdot}\pi\bbeta^{n}(\alpha^{n}_{\cdot})}_{x^{n}}
\int_{0}^{ \tau}f ( x_{t})e^{-\phi_{t} }
\,dt
$$
$$
\geq -2/n+I_{n}+v(x^{n}).
$$
The theorem is proved.
 
\mysection{Versions of Theorems \protect\ref{theorem 1.20.1}
and \protect\ref{theorem 4.18.1}
for uniformly nondegenerate case and proof of Theorem
\protect\ref{theorem 2.19.1}}  
                                              \label{section 4.14.1}
In Theorem \ref{theorem 2.19.4} below
 we suppose that 
Assumptions  \ref{assumption 1.9.1} (i) b), (ii) 
are satisfied and domain $D$
is {\em bounded\/}. We also take extensions of $\sigma,b,c$ and $f$
as in Sections \ref{section 2.25.1} and \ref{section 3.18.1}
satisfying Assumptions \ref{assumption 4.29.2}
and \ref{assumption 4.29.1}
and additionally require the extended $\sigma^{\alpha\beta}$
to also satisfy Assumption \ref{assumption 1.9.1} (iv),
perhaps with a different constant $\delta$.

Finally, we suppose that Assumptions \ref{assumption 2.23.1}
and \ref{assumption 2.23.01} are satisfied.
 
 Then take $\gamma$ and
$\lambda$ as in 
Section \ref{section 2.25.1} (and Section \ref{section 3.18.1})
and introduce
the functions $v_{\pm K}$ and $v$ as in 
Sections \ref{section 2.25.1} and \ref{section 3.18.1}.

\begin{theorem}
                                            \label{theorem 2.19.4}
We have   $v_{\pm K}\to v$ uniformly on $\bar{D}$ as $K\to \infty$.
\end{theorem} 

Proof. For $\varepsilon>0$
we construct $v_{\varepsilon,\pm K}(x)$ and $v_{\varepsilon}(x)$ 
from $\sigma,b,c^{(\varepsilon)}$,
 $f^{(\varepsilon)}$ (mollifying only the original
$c,f$ and not their extensions)
  and $g$ in the same way as $v_{\pm K}$ and $v$
were constructed from
$\sigma,b,c$, $f$,  and $g$. By Theorems \ref{theorem 1.20.1}
and \ref{theorem 4.18.1}
we have $v_{\varepsilon,\pm K}\to v_{\varepsilon}$ 
uniformly on $\bar{D}$ as $K\to \infty$
for any $\varepsilon>0$. 

Therefore, we only need to show that $|v_{\varepsilon,\pm K}-
v_{\pm K}|+|v_{\varepsilon }-
v |\leq W(\varepsilon)$, where $W(\varepsilon)$
is independent of $K$ and tends to zero as
$\varepsilon\downarrow0$.
However, by Theorem 2.2.1  of \cite{Kr77} and Lemma \ref{lemma 2.24.1} 
(see also Remarks \ref{remark 2.26.3} and \ref{remark 2.26.2})
$$
|v_{\varepsilon,\pm K}-
v_{\pm K}|+|v_{\varepsilon }-
v |\leq  N
\|\sup_{\alpha\in A,\beta\in B}|f^{\alpha\beta}-
(f^{\alpha\beta})^{(\varepsilon)}|\,\|_{L_{d}(D)}
$$
$$
+N\|\sup_{\alpha\in A,\beta\in B}|f^{\alpha\beta}
|\,\|_{L_{d}(D)}
\|\sup_{\alpha\in A,\beta\in B}|c^{\alpha\beta}-
(c^{\alpha\beta})^{(\varepsilon)}|\,\|_{L_{d}(D)}.
$$
This proves the theorem.

In the remaining part of the section   the assumption of
 Theorem  \ref{theorem 2.19.1}, that is all the assumptions
stated in Section \ref{section 2.26.3},
are supposed to be satisfied.

{\bf Proof of Theorem \ref{theorem 2.19.1}}.
For obvious reasons while proving 
the inequalities \eqref{4.2.1} and \eqref{4.2.r} in
assertions (i) and (ii) we may assume that
$g\in C^{2}(\bR^{d})$.

 (i)   
First  suppose that $D\in C^{2}$.
By Theorem 1.1 of \cite{Kr12.2}  there is a set $A_{2}$ and bounded 
continuous
functions $\sigma^{\alpha}=\sigma^{\alpha\beta}$,
$b^{\alpha}=b^{\alpha\beta}$, $c^{\alpha}=c^{\alpha\beta}$  
(independent of  $x$ and $\beta$), and $f^{\alpha\beta}\equiv0$ defined on $A_{2}$ 
such that   Assumption \ref{assumption 1.9.1} (iv)
about the uniform nondegeneracy of 
$a^{\alpha}=a^{\alpha\beta}
=(1/2)\sigma^{\alpha}(\sigma^{\alpha})^{*}$ is satisfied
for $\alpha\in  A_{2}$ (perhaps with a different constant $\delta>0$)
  and such
 that for any  $K\geq0$ the equation (the following notation 
is explained below)
\begin{equation}
                                                     \label{4.14.2}
H_{K}[u]=0
\end{equation}
(a.e.) in $D$ with boundary condition $u=g$ on $\partial D$
has a unique solution 
$$
u_{K} \in  C^{1}(\bar{D})\bigcap_{p\geq1}
  W^{2}_{p}(D) 
$$
(recall Assumption \ref{assumption 5.14.1}
 (iv) and that $g\in C^{2}(\bR^{d})$). Here              
 \begin{equation}
                                            \label{4.14.1}
H_{K}[u](x):=\max(H[u](x),
P[u](x)-K),
\end{equation}
 \begin{equation}
                                                  \label{2.8.5}
P[u](x)=\sup_{\alpha\in A_{2}}
\big[a_{ij}^{\alpha} D_{ij}u(x)
+b^{\alpha}_{i} D_{i}u(x)-c^{\alpha} u(x) \big].
\end{equation}

Observe  that
$$
\max(H[u](x),
P[u](x)-K)
$$
$$
=\max\big\{\supinf_{\alpha\in A_{1}\,\,\beta\in B}
[L^{\alpha\beta}u(x)+f^{\alpha\beta}(x)],
\supinf_{\alpha\in A_{2}\,\,\beta\in B}
[L^{\alpha\beta}u(x)+f^{\alpha\beta}(x)-K]\big\}
$$
$$
=\supinf_{\alpha\in\hat{A}\,\,\beta\in B}
\big[ L^{\alpha\beta}u( x)+f^{\alpha\beta}_{K}( x)]
\quad (f^{\alpha\beta}_{K}( x)=f^{\alpha\beta} ( x)
I_{\alpha\in A_{1}}-KI_{\alpha\in A_{2}}),
$$
where the first equality follows from the definition of $H[u]$,
\eqref{2.8.5}, and the fact that $L^{\alpha \beta}$
is independent of $\beta$ for $\alpha\in A_{2}$.

We set $u_{K}(x)=g(x)$
if $x\not\in D$.

Since $D$ is sufficiently regular by assumption, there exists
a sequence $u^{n}(x)$ of functions of class $C^{2}(\bar{D})$,
which converge to $u_{K}$ as $n\to\infty$ uniformly in $\bar{D}$ and in
$W^{2}_{p}(D)$ for any $p\geq1$. Hence, 
by Theorem 4.2 of \cite{Kr_1} we have
$$
u _{K}(x)=\infsup_{\bbeta\in\hat{\bB}\,\,\alpha\in\hat{\frA}}
E_{x}^{\alpha_{\cdot}\bbeta(\alpha_{\cdot})}\big[\int_{0}^{\tau}
 f_{K} 
( x_{t})e^{-\phi_{t}}\,dt+g(x_{\tau})e^{-\phi_{\tau}}\big].
$$
 
By Theorem \ref{theorem 2.19.4} we have that $u_{K}\to v$ uniformly on $\bar{D}$ and, since they coincide
outside $D$, the convergence is uniform on $\bR^{d}$. 
In particular,
\begin{equation}
                                                 \label{3.20.07}
v\in C(\bR^{d}).
\end{equation}

On the other hand by
assumption, \eqref{4.14.2}, and \eqref{4.14.1} we have  
$\check{H}[\psi  u_{K}]\leq 0$ (a.e. $D^{\v}$). We also know that
$u_{K}\geq v$ and, in particular, $\psi u_{K}\geq\check{v}$
on $\partial D^{\v}$ 
Furthermore, $\psi u^{n}\in C^{2}(\bar{D^{\v}})$,
$\psi u^{n}$ are $p$-insensitive by Assumption 
\ref{assumption 5.9.1},
and, for each $n$,
 the second-order derivatives of $\psi u^{n}$
are uniformly continuous in $\bar{D^{\v}})$ (because
of our assumptions on $\Pi$ and $\psi$).
Also
$\psi u^{n}$ converge to $\psi u_{K}$  as $n\to\infty$  uniformly
in $\bar{D^{\v}}$ and, as is easy to see,   
for any   $\check{x}\in D^{\v}$
$$
E^{\alpha_{\cdot}\beta_{\cdot}}_{\check{x}}
\int_{0}^{\check \tau}\big(|D^{2}_{\check{x}}( \psi   u^{n}) -
D^{2}_{\check{x}}
(  \psi   u_{K} )|
+|D_{\check{x}}(  \psi   u^{n})-D_{\check{x}}( 
 \psi  u_{K}) |\big)(
\check{x} _{t})e^{-\check{\phi}_{t}}\,dt
$$
$$
\leq NE^{\alpha_{\cdot}\beta_{\cdot}}_{\check{x}}
\int_{0}^{\check \tau}
\bigg(|D^{2}u^{n}-D^{2}u_{K}|+ |Du^{n}-Du_{K}|
+| u^{n}- u_{K}|\bigg)(\Pi(\check{x}_{t}))\,dt,
$$
where, as always, $u_{K}(\check{x})=u_{K}(\Pi(\check{x}))$
and $u^{n}(\check{x})=u^{n}(\Pi(\check{x}))$ and
 the constant $N$ depends only on  
$\|\psi, \Pi\|_{C^{1,1}}$,
$d$, and $k$. By Assumption \ref{assumption 3.25.1} (iii)  
and Theorem 2.2.1  of \cite{Kr77} the last expression tends to zero
as $n\to\infty$   uniformly with respect to $\alpha_{\cdot}\in\frA$
and $\beta_{\cdot}\in\frB$.
We also recall that the remaining parts of   
Assumption \ref{assumption 3.25.1} are imposed and
this allows us to apply  Theorem
4.1 of \cite{Kr_1}  and conclude
 that $\psi u_{K}\geq\check{v}$, which after setting
$K\to\infty$ yields $\psi v\geq\check{v}$.
Theorem
4.1 of \cite{Kr_1} also says that
$$
\psi u_{K}(\check{x})
\geq\infsup_{\bbeta\in\bB\,\,\alpha_{\cdot}\in\frA}
E_{\check{x}}^{\alpha_{\cdot}\bbeta(\alpha_{\cdot})}\big[
 \psi u_{K}(\check{x}_{\gamma\wedge\check \tau})
 e^{-\check{ \phi}_{\gamma\wedge\check \tau}-
 \psi _{\gamma\wedge\check \tau}}
$$
\begin{equation}
                                            \label{4.2.5}
+\int_{0}^{\gamma\wedge\check \tau}
\{\check{f}(p _{t},\check{x} _{t})+\lambda_{t}
 \psi u_{K}(\check{x} _{t})\}e^{-\check{ \phi}_{t}- \psi _{t}}\,dt \big].
\end{equation}
By letting $K\to\infty$ in \eqref{4.2.5}
and using the uniform convergence
of $u_{K}$ to $v$ we easily get 
the desired result  in our particular case of smooth $D$.

So far we did not use  the assumption concerning
the boundary behavior of $G$ and $\check{G}$ which we need now to 
 deal with the case of general $D$. Take an expanding sequence  
of smooth domains $D_{n}\subset D$ such that $D=\bigcup D_{n}$
and construct the functions $v^{n}$ in the same way as $v$
by replacing $D$ with $D_{n}$. We extend $v^{n}$ to $\bR^{k}$
as in \eqref{3.23.1}.

Also construct $\check{v}^{n}$ 
by replacing $D^{\v}$ with
$$
D^{\v}_{n}=D^{\v}\cap\{\check{x}:\Pi(\check{x})
\in D_{n}\}\quad(=D_{n}\quad\text{if}
\quad k =d\quad\text{and}\quad \Pi(x)\equiv x)
$$
and the boundary data   $\psi v^{n}$   
in place of $\psi v$, that is
$$
\check{v}^{n}(\check{x})
=\infsup_{\bbeta\in\bB\,\,\alpha_{\cdot}\in\frA}
E_{\check{x}}^{\alpha_{\cdot}\bbeta(\alpha_{\cdot})}
\big[\psi v^{n}(\check{x}_{\check \tau(n)} ) 
e^{-\check{\phi}_{\check \tau(n)}}
+\int_{0}^{\check \tau(n)}\check{f}(p_{t},\check{x} _{t})
e^{-\check{\phi}_{t}}\,dt\big],
$$
where 
$\check \tau^{\alpha_{\cdot}\beta_{\cdot}\check{x}}(n)$ 
is the first exit time of 
$\check{x}_{t}^{\alpha_{\cdot}\beta_{\cdot}\check{x}}$
 from $D^{\v}_{n}$.
Then by the above we have that  
\begin{equation}
                                               \label{3.20.6}
\psi v^{n}\geq \check{v}^{n}
\end{equation}
and 
$$
 \psi v^{n}(\check{x})
\geq\infsup_{\bbeta\in\bB\,\,\alpha_{\cdot}\in\frA}
E_{\check{x}}^{\alpha_{\cdot}\bbeta(\alpha_{\cdot})}\big[
 \psi v^{n}(\check{x}_{\gamma\wedge\check \tau(n)})
 e^{-\check{ \phi}_{\gamma\wedge\check \tau(n)}
- \psi _{\gamma\wedge\check \tau(n)}}
$$
\begin{equation}
                                            \label{4.2.3}
+\int_{0}^{\gamma\wedge\check \tau(n)}
\{\check{f}(p_{t},\check{x} _{t} )+\lambda_{t}
 \psi v^{n}(\check{x} _{t})\}
e^{-\check{ \phi}_{t}- \psi _{t}}\,dt \big].
\end{equation}

We now claim that, as $n\to\infty$,
\begin{equation}
                                               \label{3.20.7}
\sup_{\bR^{d}}|v^{n}-v|\to0 ,
\end{equation}
\begin{equation}
                                               \label{3.20.8}
\sup_{\bR^{k }}|\check{v}^{n}-\check{v}| \to0.
\end{equation}

 That  \eqref{3.20.7} holds is proved in \cite{Kr_1}
(see Section 6 there). 
 Owing to \eqref{3.20.7} to prove \eqref{3.20.8} it suffices to show
that uniformly in $\bR^{k }$ (notice the replacement of $v^{n}$ by $v$)
$$
\infsup_{\bbeta\in\bB\,\,\alpha_{\cdot}\in\frA}
E_{\check{x}}^{\alpha_{\cdot}\bbeta(\alpha_{\cdot})}
\big[\psi v  (\check{x}_{\check \tau(n)} ) 
e^{-\check{\phi}_{\check \tau(n)}}
+\int_{0}^{\check \tau(n)}\check{f}(p_{t},\check{x} _{t})
e^{-\check{\phi}_{t}}\,dt\big]
$$
\begin{equation}
                                              \label{3.20.9}
\to \infsup_{\bbeta\in\bB\,\,\alpha_{\cdot}\in\frA}
E_{\check{x}}^{\alpha_{\cdot}\bbeta(\alpha_{\cdot})}
\big[\psi v (\check{x}_{\check \tau } ) 
e^{-\check{\phi}_{\check \tau }}
+\int_{0}^{\check \tau }\check{f}(p_{t},\check{x} _{t})
e^{-\check{\phi}_{t}}\,dt\big]
\end{equation}
(recall that $\check \tau^{\alpha_{\cdot}\beta_{\cdot}\check{x}}$ 
is the first exit time of 
$\check{x}_{t}^{\alpha_{\cdot}\beta_{\cdot}\check{x}}$ 
from $D^{\v}$). Both sides of \eqref{3.20.9} coincide
if $\check x\not\in D^{\v}$. Therefore, we need to prove the
uniform convergence only in $D^{\v}$.

Here $v\in C(\bar{D})$ and it is convenient to prove
\eqref{3.20.9} just for any such $v$, regardless of its
particular construction. In that case, 
relying on Assumption \ref{assumption 5.9.1},
as in the proof
of Theorem 2.2 of \cite{Kr_1} (see Section 6
there), we reduce
our problem to proving that uniformly in $D^{\v}$ 
$$
\hat{v}_{n}(x):=\infsup_{\bbeta\in\bB\,\,\alpha_{\cdot}\in\frA}
E_{\check{x}}^{\alpha_{\cdot}\bbeta(\alpha_{\cdot})}
\int_{0}^{\check \tau(n)}\check{f}(p_{t},\check{x} _{t})
e^{-\check{\phi}_{t}}\,dt  
$$
$$
\to \hat{v}(x):=\infsup_{\bbeta\in\bB\,\,\alpha_{\cdot}\in\frA}
E_{\check{x}}^{\alpha_{\cdot}\bbeta(\alpha_{\cdot})}
 \int_{0}^{\check \tau }\check{f}(p_{t},\check{x} _{t})
e^{-\check{\phi}_{t}}\,dt
$$
with perhaps modified $\check{f}$ still satisfying 
\eqref{6.14.4} and \eqref{3.27.04} and satisfying
 Assumptions \ref{assumption 3.25.1} (i), (ii)
with (modified $\bar{f}_{\varepsilon}$).  

For $\varepsilon>0$ introduce
$$
N_{\varepsilon}=\sup_{(\alpha,\beta,\check x)\in
A\times B\times D^{\v}}|
\bar f_{\varepsilon}^{\alpha\beta}(\check x)|
$$
and observe that
$$
|\check v(\check x)-\check v_{n}(\check x)|
\leq \delta_{1}^{-1}I_{n}(x),
$$
where
$$
\delta_{1}I_{n}(x)=\delta\sup_{\alpha_{\cdot}\in\frA,
\beta_{\cdot}\in\frB}E^{\alpha_{\cdot}\beta_{\cdot}}_{\check x}
\int_{\check \tau_{n}}^{\check\tau}
|f(p_{t},\check x_{t})|e^{-\check \phi_{t}}\,dt
$$ 
$$
\leq 
\sup_{\alpha_{\cdot}\in\frA,
\beta_{\cdot}\in\frB}E^{\alpha_{\cdot}\beta_{\cdot}}_{\check x}
\int_{\check \tau_{n}}^{\check\tau}
|\bar f  (\check x_{t})|e^{-\check \phi_{t}}\,dt
\leq N_{\varepsilon}
\sup_{\alpha_{\cdot}\in\frA,
\beta_{\cdot}\in\frB}E^{\alpha_{\cdot}\beta_{\cdot}}_{\check x}
\int_{\check \tau_{n}}^{\check\tau}
e^{-\check \phi_{t}}\,dt+J_{n}(x),
$$
where
$$
J_{n}(x)
=\sup_{\alpha_{\cdot}\in\frA,
\beta_{\cdot}\in\frB}E^{\alpha_{\cdot}\beta_{\cdot}}_{\check x}
\int_{0}^{\check\tau}
|\bar f  (\check x_{t})-\bar f_{\varepsilon}
(\check x_{t})|e^{-\check \phi_{t}}\,dt. 
$$

By Assumption \ref{assumption 3.25.1} (ii) we have that
$J_{n}(x)\to0$ as $\varepsilon\downarrow0$ uniformly
in $D^{\v}$ (this is the only place where we use
the uniformity in \eqref{3.27.2}). Furthermore, by   Lemma 5.1 of \cite{Kr_1}  
\begin{equation}
                                               \label{5.29.1}
\sup_{\alpha_{\cdot}\in\frA,
\beta_{\cdot}\in\frB}E^{\alpha_{\cdot}\beta_{\cdot}}_{\check x}
\int_{\check \tau_{n}}^{\check\tau}
e^{-\check \phi_{t}}\,dt\leq\sup_{\partial D^{\v}_{n} }\check G.
\end{equation}
As is easy to check $\Pi(\partial D^{\v}_{n})
\subset \partial D_{n}$, so that, if we have a sequence of points
$\check x_{n}\in \partial D^{\v}_{n}$, then
$\dist(\Pi(\check x_{n}),\partial D)\to0$ as 
$n\to\infty$. It follows 
by  Assumption \ref{assumption 5.24.1} that the right-hand side
of \eqref{5.29.1} goes to zero as $n\to\infty$.
 This proves that $I_{n}(x)\to0$ uniformly 
in $D^{\v}$, yields \eqref{3.20.9} and
\eqref{3.20.8} and along with \eqref{3.20.7}
and \eqref{3.20.6} proves that $\psi v\geq\check{v}$.
One passes to the limit in \eqref{4.2.3} similarly
and this
finally brings the proof of assertion (i)
to an end.

(ii) 
As above first  suppose that $D\in C^{2}$. By
Theorem 1.3 of \cite{Kr12.2}   there is a set $B_{2}$ and bounded 
continuous
functions $\sigma^{\beta}=\sigma^{\alpha\beta}$,
$b^{\beta}=b^{\alpha\beta}$, $c^{\beta}=c^{\alpha\beta}$  
(independent of  $x$ and $\alpha$), and 
$f^{\alpha\beta}\equiv0$ defined on $B_{2}$ 
such that   Assumption \ref{assumption 1.9.1} (iv)
about the uniform nondegeneracy of 
$a^{\beta}=a^{\alpha\beta}
=(1/2)\sigma^{\beta}(\sigma^{\beta})^{*}$ is satisfied
for $\beta\in  B_{2}$ (perhaps with a different constant $\delta>0$)
  and such
 that for any  $K\geq0$ the equation (the following notation 
is explained below)
$$
H_{-K}[u]=0
$$
(a.e.) in $D$ with boundary condition $u=g$ on $\partial D$
has a unique solution 
$$
u_{-K} \in  C^{1}(\bar{D})\bigcap_{p\geq1}
  W^{2}_{p}(D) .
$$
 Here              
 $$
H_{-K}[u](x):=\max(H[u](x),
P[u](x)+K),
$$
$$
P[u](x)=\inf_{\beta\in B_{2}}
\big[a_{ij}^{\beta}D_{ij}u(x)
+b^{\beta}_{i}D_{i}u(x)-c^{\beta}u(x)\big].
$$

We introduce
$$
f^{\alpha\beta}_{K}( x)=f^{\alpha\beta}( x)
I_{\beta\in B_{1}}+KI_{\beta\in B_{2}}.
$$
and note that
$$
H_{-K}[u](x):=\supinf_{\alpha\in A\,\,\beta\in \hat{B}}
\big[ L^{\alpha\beta}u( x)+f^{\alpha\beta}_{K}( x)]
$$
$$
=\sup_{\alpha\in A }\min\big\{\inf_{ \beta\in B}
[  L^{\alpha\beta}u(x)+ f^{\alpha\beta}(x)],
 \inf_{\beta\in B_{2}}
[L^{\alpha\beta}u(x)+f^{\alpha\beta}(x)+K]\big\}
$$
$$
=\min(H[u](x),
P[u](x)+K).
$$
After that it suffices to repeat the above proof
relying   
 again on Theorem \ref{theorem 2.19.4}.

The theorem is proved.

\mysection{Proof  of Theorems  \protect\ref{theorem 1.14.1} 
and \protect\ref{theorem 2.18.1}}

                                          \label{section 4.16.2}

Here all assumptions of Theorem \ref{theorem 1.14.1}
are supposed to be satisfied.

{\bf Proof of Theorem \ref{theorem 1.14.1}}.
  If the functions $(c,f)^{\alpha\beta}(x)$ 
are bounded on $A\times
B\times\bR^{d}$, then according to Remark \ref{remark 4.2.1}
assertion (ii) of Theorem \ref{theorem 1.14.1}
follows immediately from Theorem \ref{theorem 2.19.1}. The continuity
of $v$ also follows from the proof of Theorem \ref{theorem 2.19.1}.

In the general case, for $\varepsilon>0$,
define
$$
(c_{\varepsilon},f_{\varepsilon})^{\alpha\beta}(x)=
 (  c^{\alpha\beta},
  f^{\alpha\beta})^{(\varepsilon)}
(x) 
$$
and construct $v_{\varepsilon}(x)$ from $\sigma,b,c_{\varepsilon }$,
 $f_{ \varepsilon }$, and $g$ in the same way as $v$ was constructed from
$\sigma,b,c$, $f$,  and $g$. By the above  \eqref{1.14.1}
holds if we replace $f$ and $c$ with $f_{\varepsilon }$
and $c_{\varepsilon }$ respectively. 

We first take $\lambda\equiv0 $ and $\gamma^{\alpha_{\cdot}\beta_{\cdot}}
=\tau^{\alpha_{\cdot}\beta_{\cdot}x}$ in the counterpart
of \eqref{1.14.1} corresponding to $v_{\varepsilon}$.
Then by  Theorem 2.2.1  of \cite{Kr77} and Lemma \ref{lemma 2.24.1} 
(see also Remarks \ref{remark 2.26.3} and \ref{remark 2.26.2})
$$
|v(x)-v_{\varepsilon}(x)|\leq  N
\|\sup_{\alpha\in A,\beta\in B}|  f^{\alpha\beta}-
( f^{\alpha\beta})^{(\varepsilon)}|\,\|_{L_{d}(D)}
$$
$$
+N\|\sup_{\alpha\in A,\beta\in B}| f^{\alpha\beta}
|\,\|_{L_{d}(D)}
\|\sup_{\alpha\in A,\beta\in B}|  c^{\alpha\beta}-
( c^{\alpha\beta})^{(\varepsilon)}|\,\|_{L_{d}(D)}.
$$
It follows by Assumption \ref{assumption 1.9.1} (iii) that
$v_{\varepsilon}\to v$ uniformly on $\bar{D}$ and $v$
is continuous in $\bar{D}$.
After that we easily pass to the limit in the counterpart
of \eqref{1.14.1} corresponding to $v_{\varepsilon}$
for arbitrary
$\lambda $ and $\gamma$ again on the basis of Lemma
\ref{lemma 2.24.1}. The theorem is proved.  

{\bf Proof of Theorem \ref{theorem 2.18.1}}.
We know from \cite{Kr12.2} (see Remark 1.3 there)
 that $u_{K}$ introduced in the proof
of Theorem \ref{theorem 2.19.1} (see Section \ref{section 4.14.1})
satisfies an elliptic equation 
$$
a^{K}_{ij}D_{ij}u_{K}+b^{K}_{i}D_{i}u_{K}-c^{K}u_{K}+f^{K}=0,
$$ 
where $(a^{K}_{ij})$
satisfies the uniform nondegeneracy condition
(see Assumption \ref{assumption 1.9.1} (iv)) with a constant
$\delta_{1}=\delta_{1}(\delta,d)>0$, $|b^{K}|,c^{K}$  are bounded
by a constant depending only on $K_{0}$, $\delta$, and $d$,   
$ c^{K}\geq0$ and
 $$
|f^{K}|\leq \sup_{\alpha,\beta}| f^{\alpha\beta}|.
$$
Then according to classical results (see, for instance,
\cite{GT} or \cite{Kr85}) there exists a constant $\theta\in(0,1)$
depending only on $\delta_{1}$ and $d$, that is on $\delta$ and $d$,
such that
for any subdomain $D'\subset \bar{D'}
\subset D$ and $x,y\in D'$
we have
\begin{equation}
                                                      \label{4.14.4}
|u_{K}(x)-u_{K}(y)|\leq N|x-y|^{\theta},
\end{equation}
where $N$ depends only on $\delta$, $d$,
the distance between the boundaries
of $D'$ and $D$, on the diameter of $D$, and on $K_{0}$.
It is seen that \eqref{4.14.4} will be preserved as we let $K
\to\infty$ and then perform all other steps in the above proof
of Theorem \ref{theorem 1.14.1} which will lead us to the desired
result. The theorem is proved.


\begin{thebibliography}{mm}

 
\bibitem{BL08} R. Buckdahn and J. Li, {\em
Stochastic differential games and viscosity solutions of 
Hamilton-Jacobi-Bellman-Isaacs equations\/}, SIAM J. Control Optim.,
Vol. 47 (2008), No. 1, 444--475.


\bibitem{FS89} W. H.
Fleming and P. E. Souganidis, {\em On the existence of
value functions of two-player, zero-sum stochastic differential games\/},
Indiana Univ. Math. J., Vol. 38 (1989), No. 2, 293--314. 
 
 
 \bibitem{GT} D. Gilbarg and N.S. Trudinger, ``Elliptic Partial 
Differential Equations of Second Order",
Series: Classics in Mathematics,  Springer,
Berlin, Heidelberg, New York,  2001.

\bibitem{Ko09} J. Kovats, {\em
Value functions and the Dirichlet problem for Isaacs equation in a smooth
domain\/}, Trans. Amer. Math. Soc., Vol. 361 (2009), No. 8,
4045--4076. 

 \bibitem{Kr70}  N.V. Krylov,   {\em
    On a problem with two free boundaries for an elliptic
   equation and optimal stopping of a Markov process\/},
Doklady Academii Nauk SSSR,
 Vol. 194 (1970), No. 6, 1263--1265 in Russian; English
translation:  Soviet Math. Dokl.,
   Vol. 11
(1970),
No. 5, 1370--1372.

\bibitem{Kr77}   N.V. Krylov, ``Controlled diffusion processes'', 
Nauka, Moscow,  1977 in Russian; English translation  Springer,
1980.
 

\bibitem{Kr85} N.V. Krylov,
``Nonlinear elliptic and parabolic equations of second
  order'',  Nauka, Moscow, 1985 in Russian; English translation
   Reidel, Dordrecht, 1987.

 
\bibitem{Kr95}  N.V. Krylov,
``Introduction to the theory of diffusion processes'',   Amer.
Math. Soc., Providence, RI, 1995.

  

\bibitem{Kr12.2}  N.V. Krylov,
{\em On the existence of smooth 
solutions for fully nonlinear elliptic
equations with measurable ``coefficients''
without convexity assumptions\/},
http://arxiv.org/abs/1203.1298

\bibitem{Kr_1}  N.V. Krylov,
{\em On the dynamic programming principle for
uniformly nondegenerate stochastic differential
games in domains\/}, http://arxiv.org/abs/1205.0048

\bibitem{Ni88} M. Nisio, {\em
Stochastic differential games and viscosity 
solutions of Isaacs equations\/}
Nagoya Math. J., Vol. 110 (1988), 163--184. 

 \bibitem{Sw96} A. \'Swi{\c e}ch, {\em Another approach to the existence of value
functions of stochastic differential games\/}, J. Math.
Anal. Appl., Vol. 204 (1996), No. 3, 884--897.
\end{thebibliography}
\end{document}